\documentclass[11pt]{article}
\usepackage{latexsym,amsmath,geometry, graphicx}
\usepackage{amssymb}
\usepackage{ upgreek, color}

\newcommand \nc{\newcommand}
\newtheorem{theorem}{Theorem}[section]
\newtheorem{lemma}[theorem]{Lemma}

\newtheorem{remark}[theorem]{Remark}

\nc{\ba}{\begin{array}}\nc{\ea}{\end{array}}
\nc{\be}{\begin{eqnarray}}\nc{\ee}{\end{eqnarray}}
\nc{\beq}{\begin{equation}}\nc{\eeq}{\end{equation}}
\nc{\bex}{\begin{eqnarray*}}\nc{\eex}{\end{eqnarray*}}
\nc{\btm}{\begin{theorem}} \nc{\etm}{\end{theorem}}
\nc{\blm}{\begin{lemma}} \nc{\elm}{\end{lemma}}
\nc{\R}{\mathbb{R}} \nc{\va}{\varepsilon} \nc{\ls}{\limits}

\def\e{\varepsilon}

\def\nn{\nonumber}
\def\pf{\noindent{\bf Proof.\quad}}\def\endpf{\hfill$\Box$}

\begin{document}
\title{Finite time singularities for hyperbolic systems}
\author{Geng Chen
\thanks{School of Mathematics,
Georgia Institute of Technology, Atlanta, GA 30332 USA, gchen73@math.gatech.edu}
\qquad
Tao Huang
\thanks{Department of Mathematics,
Pennsylvania State
University, University Park, PA 16802, USA, txh35@psu.edu}
\qquad
Chun Liu
\thanks{Department of Mathematics,
Pennsylvania State
University, University Park, PA 16802, USA, cxl41@psu.edu}
}
\date{}
\maketitle

\abstract{
In this paper, we study the formation of finite time singularities in the form of super norm blowup
for a spatially inhomogeneous hyperbolic system. The system is related to the variational wave equations as those in \cite{ghz}.
The system posses a unique $C^1$ solution before the emergence of 
vacuum in finite time, for given initial data that are smooth enough, bounded and uniformly away from vacuum. 
At the occurrence of blowup, the density becomes zero, while the momentum stays finite, however 
the velocity and the energy are both infinity. 
}

\section{Introduction}
\setcounter{equation}{0}
\setcounter{theorem}{0}
In this paper, we consider the following Cauchy problem of spatially inhomogeneous hyperbolic partial differential equations:
\beq
\label{LCEuler}
\left\{
\begin{array}{l}
\rho_t+(\rho u)_x=0\\
(\rho u)_t+\bigl(\rho u^2- c^2(x)\rho^{-1}\bigr)_x=-c(x) c'(x) \rho^{-1}\\
(\rho, u)|_{t=0}=(\rho_0, u_0),
\end{array}
\right.
\eeq
where $\rho(x,t):\mathbb R\times[0,+\infty)\rightarrow \mathbb [0,\infty)$ is the density, $u(x,t):\mathbb R\times[0,+\infty)\rightarrow \mathbb R$ is the velocity, $\rho_0, u_0$ are given initial data that will be specified later and $c(x):\mathbb R\rightarrow \mathbb R^+$ is a given function satisfying
\beq\label{c_con}
c(x)\in C^{2}, \quad 0<c_0\leq c(x)\leq C_0<+\infty, \quad |c'(x)|\leq C_1<+\infty,\quad c'(x)\not\equiv0,
\eeq
for some constants $c_0$, $C_0$ and $C_1$.
It is easy to verify that the smooth solutions of (\ref{LCEuler}) satisfy an
energy conservation law
\beq\label{energy}
E_t +q_x=0
\eeq
with specific energy (entropy)
\beq\label{energy_intr}
E=\frac{1}{2}\rho u^2+\frac{1}{2}c^2 \rho^{-1}
\eeq
and entropy flux
\beq\label{flux_intr}
q=\frac{1}{2}u^3\rho-\frac{1}{2}c^2\rho^{-1}u\,.
\eeq

\subsection{Inhomogeneous linear elasticity}
System \eqref{LCEuler} can be viewed as an Eulerian description of inhomogeneous linear elasticity.

\begin{figure}[htb]
\centering
\includegraphics[width=16cm,height=4cm]{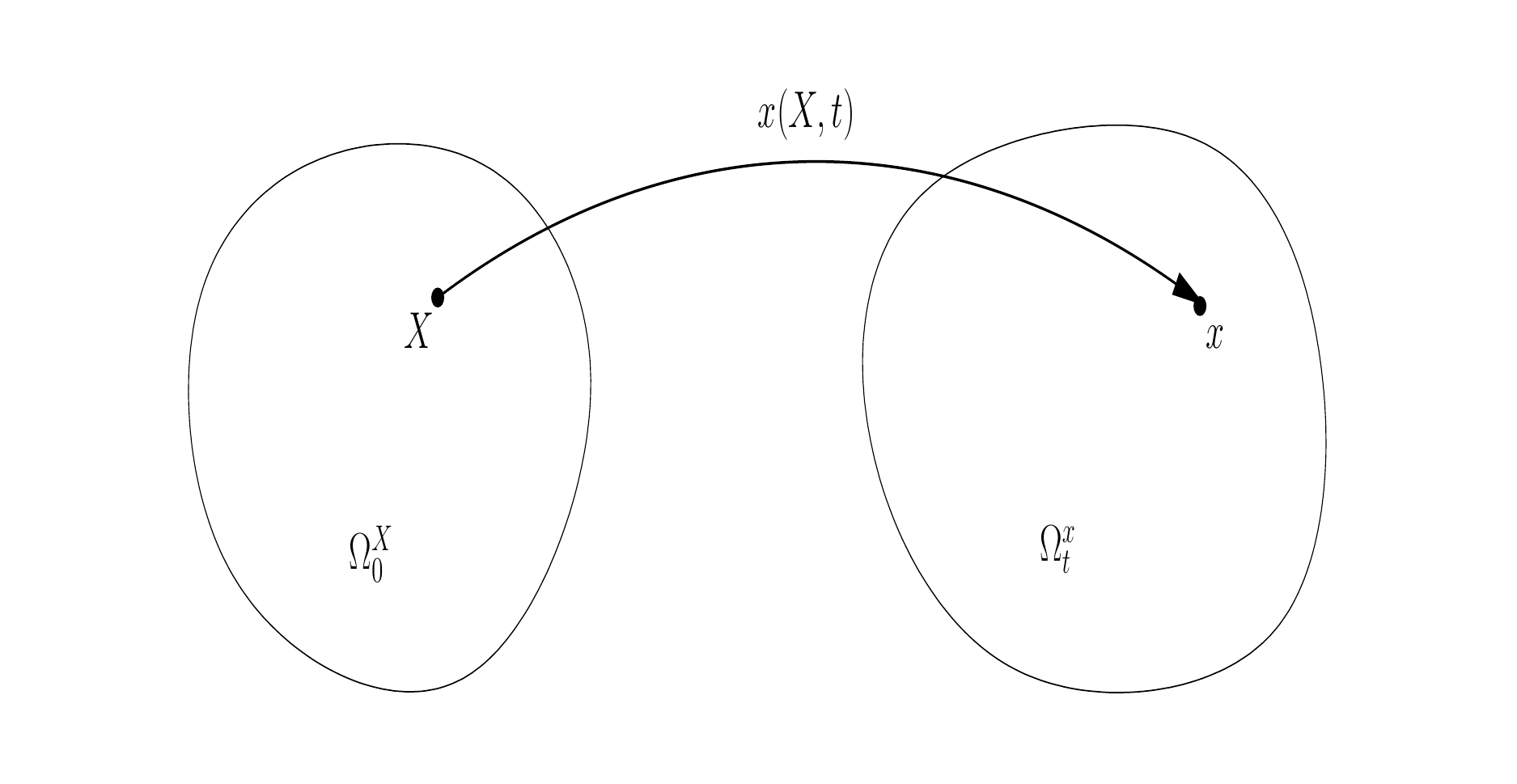}
\caption{A flow map}
\label{flowmap}
\end{figure}

 \bigskip
\noindent{\bf Basic mechanics.} For any smooth domain $\Omega_0^{X}\subset\mathbb R^n$ with $n=1, 2$, or $3$, 
let $X\in \Omega_0^{X}$ denotes the Lagrangian coordinates. The flow map $x(X,t):\Omega_0^X\rightarrow\Omega_t^x$ in Figure \ref{flowmap} satisfies (see \cite{arnold} for details)
\begin{equation}
\label{sn3}
\left\{
\begin{array}{rl}
\displaystyle \frac{dx}{dt}&=u(x(X,t),t)\\
x(X,0)&=X\,.
\end{array}
\right.
\end{equation}

\noindent Furthermore, let
\begin{equation}\label{FLang}
\displaystyle \tilde{F}(X,t)=\frac{\partial x(X,t)}{\partial X}
\end{equation}
be the deformation matrix associate with the flow map and define the Eulerian quantity (push forward)
\begin{equation}\label{FEuler}
F(x(X,t),t)=\tilde{F}(X,t).
\end{equation}
In case of $F$ preserving the sign, e.g. $\det\,\tilde F>0$, not only $\tilde F$ is an invertible matrix, $x(X,t)$ preserves the orientation. By \eqref{FLang}, \eqref{FEuler} and direct calculation, we obtain that $F$ satisfies the following kinematic relation (chain rule)  (c.f. \cite{lin-liu-zhang})
\begin{equation}
\label{eqnF}
F_t+u\cdot\nabla_x F=\nabla_x u\cdot F.
\end{equation}
Let $\rho(x,t):\Omega_t^x\times[0,+\infty)\rightarrow \mathbb R^+$ be the density of mass with initial data $\rho_{0}(X)=\rho(x(X,0),0)$. The usual conservation of mass equation
\begin{equation}
\label{eqnconsermass}
\rho_t+\nabla_x\cdot(\rho u)=0
\end{equation}
is equivalent to
\begin{equation}\label{equationRho}
\rho(x(X,t),t)=\frac{\rho_0(X)}{\det\,\tilde F(X,t)}.
\end{equation}


 \bigskip
\noindent{\bf Energetic variational approaches.} By the first and second laws of thermodynamics, one can start with energy 
law for a conservative system:
\begin{equation}
\label {genenergylaw}
\frac{d}{dt}\,\int_{\Omega_0^X}\mathcal K(X,t,x_t)+\mathcal W(X,t,x,\tilde F)\,dX=0,
\end{equation}
where $\mathcal K(X,t,x_t)$ denotes the kinetic energy and $\mathcal W(X,t,x,\tilde F)$ denotes the free energy. 
Specifically, one can consider the following simple forms of kinetic energy and internal energy (for inhomogeneous linear elasticity) 
\begin{equation}
\label {sn2}
\frac{d}{dt}\frac{1}{2}\,\int_{\Omega_0^X}\rho_{0}(X)|x_t(X,t)|^2+c^2(x)|\tilde F(X,t)|^2\,dX=0,
\end{equation}
where $c(x):\Omega_t^x\rightarrow\mathbb R^+$ is a given scalar function.
The energy law \eqref{sn2} has been widely used to describe the elasticity in an inhomogeneous medium, that includes the coupling and competition between the kinetic energy and (linear) elastic energy (c.f \cite{liu-walkington}, \cite{lin-liu-zhang} and \cite{lei-liu-zhou}). In particular, when the space dimension $n=1$ and the initial data $\rho_0(X)=1$, we have
\begin{equation}
\label{sn10}
\rho(x(X,t),t)=\frac{\rho_0(X)}{\det\tilde F(X,t)}=\frac{1}{\tilde F(X,t)}>0.
\end{equation}
Therefore, to certain degree, in one space dimension one cannot distinguish elastic energy which depends only on $\tilde F$ and the conventional free energy of fluid which is a function of $\rho$ (c.f. \cite{Dafermos}).  

System \eqref{LCEuler} can be derived from the energy law \eqref{sn2} by using the energetic variational methods under Eulerian coordinates (c.f. \cite{huang-liu-weber}). For completeness of our paper, we sketch main steps of derivation as those in \cite{huang-liu-weber}. 

For any $T>0$, the action of the system is
$$A(x)=\frac12\int_{0}^{T}\int_{\Omega}\rho
_0(X)|x_t(X,t)|^2-c^2(x)|\tilde F(X,t)|^2\,dXdt. $$
By the least action law, the variation of $A(x)$ with respect to $x$ under Lagrangian coordinates can be calculated as follows:
\begin{equation}
\label{sn4}
\begin{split}
0&=\left.\frac{d}{d\epsilon}\right|_{\epsilon=0}A(x+\epsilon y)\\
&=\int_{0}^{T}\int_{\Omega_0^X}\rho
_0(X)x_t\cdot y_t-c(x)\,|\tilde F|^2\,(y\cdot\nabla_x) c(x)-c^2(x)\tilde F:\frac{\partial y}{\partial X}\,dXdt\\
&=-\int_{0}^{T}\int_{\Omega_0^X}\rho
_0(X)x_{tt}\cdot y+c(x)\,|\tilde F|^2\,(y\cdot\nabla_x) c(x)-\mbox{div}_{X}\left(c^2(x)\tilde F\right)y\,dXdt,
\end{split}
\end{equation}
for any $y(X,t)\in C^{\infty}_0(\Omega_0^X\times(0,T))$. Here $A:B=A_{ij}B_{ij}$ denotes the inner product of two matrices. 
Therefore, the Euler-Lagrange equation of energy law \eqref{sn2} under Lagrangian coordinates for any $n=1,2$ or $3$, is
\begin{equation}
\label{sn5}
\begin{split}
\rho_0(X)x^i_{tt}+c(x)\tilde F^2\nabla_{x_i}c(x)-\nabla_{X_j}(c^2(x)\tilde F_{ij})=0.
\end{split}
\end{equation}

\begin{remark}\label{int-rm1dLan}
When the space dimension $n=1$ and the initial data $\rho_0(X)=1$, equation \eqref{sn5} can be written as nonlinear wave equation
\begin{equation}
\label{sn6}
\begin{split}
x_{tt}(X,t)-c(x)(c(x)\tilde F(X,t))_{X}=0,
\end{split}
\end{equation}
which is exactly a special case of one dimensional variational wave equation modeling nematic liquid crystal dynamics. We will provide more details later. 
\end{remark}

Let $\tilde{y}(x(X,t),t)=y(X,t)$ be the pull-forward quantity of $y$ from the Lagrangian to Eulerian coordinates. Integrating by parts with respect to $t$, changing variables and integrating by parts with respect to $x$, one can obtain 
\begin{equation}
\label{sn7}
\begin{split}
0&=\left.\frac{d}{d\epsilon}\right|_{\epsilon=0}A(x+\epsilon y)\\
&=-\int_{0}^{T}\int_{\Omega_t}\rho(x,t)(u_t+u\cdot\nabla_x u)\tilde{y}+\frac{c(x)|F|^2(\tilde{y}\cdot\nabla_x) c(x)}{\det F}-\mbox{div}_{x}\left(\frac{c^2(x)FF^T}{\det F}\right)\tilde{y}   \,dxdt,
\end{split}
\end{equation}
where \eqref{sn3} and \eqref{equationRho} have been used in changing variables.
Therefore, the Euler-Lagrange equation of energy law \eqref{sn2} in the Eulerian coordinates is
\begin{equation}
\label{sn8}
\begin{split}
\rho(x,t)(u_t+u\cdot\nabla_x u)+\frac{c(x)|F|^2\nabla_x c(x)}{\det F}-\mbox{div}_{x}\left(\frac{c^2(x)FF^T}{\det F}\right)=0.
\end{split}
\end{equation}
When the space dimension $n=1$, combining \eqref{sn8}, \eqref{eqnconsermass} and \eqref{eqnF}, we have the following coupled dynamic system 
\begin{equation}
\label{sn9}
\begin{split}
\rho_t+(\rho u)_x&=0\\
F_t+uF_x&=u_xF\\
\rho(u_t+uu_x)&
=c(x)\left(c(x)F\right)_x.
\end{split}
\end{equation}
%
%
\begin{remark}
\label{int-rm1deuler}
When the space dimension $n=1$, $\Omega_t^x=\mathbb R$, the initial data $\rho_0(X)=1$ $\left(i.e.\ F=\frac{1}{\rho}\right)$, system \eqref{sn9} becomes
\begin{equation}
\label{sn11}
\begin{split}
\rho_t+(\rho u)_x&=0\\
\rho(u_t+uu_x)&=c(x)\left(\frac{c(x)}{\rho}\right)_x,
\end{split}
\end{equation}
which is exactly \eqref{LCEuler}. And the corresponding energy law \eqref{sn2} becomes
\begin{equation}
\label {1denergy}
\frac{d}{dt}\frac{1}{2}\,\int_{\mathbb R}\rho(x,t)|u(x,t)|^2+c^2(x)\rho^{-1}(x,t)\,dx=0.
\end{equation}
\end{remark}


\bigskip
From Remark \ref{int-rm1dLan} and Remark \ref{int-rm1deuler}, when the space dimension $n=1$, the initial data $\rho_0(X)=1$, $\det F>0$ and the solutions for both systems are smooth enough, \eqref{sn11} and variational wave equation \eqref{sn6} are formally equivalent systems under different coordinates. However, in general when one looks at \eqref{sn11} and \eqref{sn6} in weak form, they can be different since the deformation matrix $F$ is not always invertible (when singularities occur). 

The equation \eqref{sn6} is a special case of variational wave equation modeling nematic liquid crystal. In \cite{AH} such variational wave equation was first investigated in any dimensions when people were trying to find the minimal of the following energy
\begin{equation}
\label{sn-in1}
\begin{split}
\int_{\Omega}|\mathbf n_t|^2-W(\mathbf n,\nabla \mathbf n)\,dy=0,
\end{split}
\end{equation}
where $|\mathbf n|=1$ and 
$$
W(\mathbf n,\nabla \mathbf n)=\alpha|\mathbf n\times(\nabla\times\mathbf n)|^2+\beta(\nabla \cdot\mathbf n)^2+\gamma(\mathbf n\cdot\nabla\times\mathbf n)^2+ \eta \bigl[    {\text{tr}}(\nabla {\mathbf n})^2- 
(\nabla\cdot {\mathbf n})^2 \bigr]\,.
$$
Here $\alpha$, $\beta$, $\gamma$ and $\eta$ are all positive viscosity constants,
and $W$ is the Oseen-Frank potential for nematic liquid crystal (c.f. \cite{AH}).
When $\mathbf n$ only depends on a single space variable $X$ and  
$$
\mathbf n=\cos \phi(X,t)\mathbf e_{y_1}+\sin\phi(X,t)\mathbf e_{y_2}\qquad \text{(planar\ deformation)},
$$
where $\mathbf e_{y_1}$ and $\mathbf e_{y_2}$ are the coordinate vectors in the $y_1$ and $y_2$ directions, respectively.
The Euler-Lagrange equation of  \eqref{sn-in1} was given in \cite{AH} as follows 
\begin{equation}
\label{equationphi}
\phi_{tt}-c(\phi)(c(\phi)\phi_X)_X=0,
\end{equation}
with
$$
c^2(\phi)=\alpha\cos^2(\phi)+\beta\sin^2(\phi).
$$
It is obvious that \eqref{equationphi} is exactly \eqref{sn6} with $\phi$ replaced by $x$.
In \cite{BZ},  Bressan and Zheng have established the global existence of energy conservative weak solutions for (\ref{sn6})) by introducing new energy-dependent coordinates (see also \cite{HR}). The solutions are locally H\"older continuous with exponent $\frac{1}{2}$. For general $W(\bf{n}, \nabla \bf{n})$ in one space dimension, the existence of weak solutions has been studied by a series of papers \cite{ZZ10,ZZ11,CZZ12}.

We really need to point out that the singularity formation for (\ref{sn6}) has been first studied by Glassey, Hunter and Zheng 
in their seminal work \cite{ghz}, in which a gradient blowup example has been provided. When there is a damping term in \eqref{sn6}, a similar gradient blowup example is provided in \cite{GCZ}. In \cite{GCZ,ghz}, the singularities they construct are "kink" solutions instead of shock waves constructed for systems of conservation laws including at least one genuinely nonlinear characteristic family \cite{lax, john, G3, G5, G6, G8}. We will provide more details in Remark \ref{int-rmelast}, Remark \ref{remark1} and Remark \ref{remarkRho0}.

Our first main result is for  \eqref{LCEuler}  with initial data $\rho_0=1$ describing the inhomogeneous elastic flow.

\begin{theorem}
\label{main-elast}
There exists a function $c(x)\in C^2$ (given in \eqref{elas1_c0}) satisfying \eqref{c_con} and a finite time $T=O(1)\e^{-2}>0$, such that the Cauchy problem of \eqref{LCEuler} with initial data $(\rho_0, u_0)=(1,c(x))$ has a unique $C^1$ solution $(\rho(x,t), u(x,t))$ on $(x,t)\in\mathbb{R}\times [0,T)$. Moreover, there exists a point 
$(x^*, T)$ with $|x^*|<O(1)\e^{-1}$ at which the solution satisfies
\beq
\lim_{t\rightarrow{T^-}}u(x^*,t)=\infty \quad \lim_{t\rightarrow{T^-}}\rho(x^*,t)=0 \quad  \lim_{t\rightarrow{T^-}}\rho u(x^*,t)=B\quad \text{and}\quad
\lim_{t\rightarrow{T^-}}E(x^*,t)=\infty,
\eeq
where $B$ is a finite constant. 
The solutions $\rho(x,t)$ and $-u(x,t)$ have uniform upper bounds 
on $(x,t)\in\mathbb{R}\times [0,T)$.
\end{theorem}

\begin{remark}
\label{int-rmelast}
We have several remarks for Theorem \ref{main-elast}:
\begin{itemize}
\item[(1).] In our example, the characteristic speeds of two families are uniformly away from each other,
or in another word, system \eqref{LCEuler} is uniformly strictly hyperbolic (see Lemma \ref{elas_lemmaR} and Lemma \ref{elas_lemmaS}). 

\item[(2).] In our example, the $L^{\infty}$ blowup happens at the same time of $C^1$ blowup,
since two characteristic families for \eqref{LCEuler} are both linearly degenerate in the definition of Lax when $c(x)$ is constant. This is a {\em rare} case for systems of hyperbolic conservation laws.
In fact, for \eqref{LCEuler-1}, this happens only when $\gamma=-1$.

\item[(3).] By \eqref{equationRho}, \eqref{sn6} and the last remark,
one could see that the $C^1$ blowup of the flow map essentially indicates the vacuum formation for \eqref{LCEuler} when the initial data $\rho_0=1$ and $F>0$, through the transformation from Lagrangian coordinates to Eulerian coordinates. 

The construction of our example is motivated by the pioneer work by Glassey, Hunter and Zheng \cite{ghz} on the singularities of variational wave equation \eqref{sn6} (with $F=x_X$) in one dimension. However, the blowup example constructed in  \eqref{sn6} does not satisfy the restrictions $\rho_0=1$ (i.e. $x_X=1$ initially) and $F=x_X>0$, hence the system cannot be transformed to \eqref{LCEuler} in this situation. By introducing new techniques, in Theorem \ref{main-elast}, we construct a blowup example satisfying all these restrictions. We do adopt important ideas from \cite{ghz}, while the restrictions for inhomogeneous elastic flow make the construction much more complicated than the example in \cite{ghz}.

\item[(4).] For any positive constant $K$, the total energy of the solution for any time until the blowup when $x\in[-K,K]$ is bounded by a constant depending on $K$, although at the time of blowup, the energy concentrates, i.e. energy density is infinity, somewhere.

\item[(5).] It is an interesting question to consider more general assumptions on $c$. The numerical experiments in \cite{huang-liu-weber} have indicated such blowup might happen for more general cases.

\item [(6).] When the initial data $\rho_0\not\equiv1$, one needs to investigate the system \eqref{sn9} instead of \eqref{LCEuler}. It is also an interesting and challenging question.

\end{itemize}

\end{remark}

\subsection{ Isentropic duct flow for Chaplygin gas dynamics}
System \eqref{LCEuler} has applications in various fields. We can reformulate the spatially inhomogeneous system \eqref{LCEuler} by setting
\beq\label{rho}
\rho=c(x)\bar\rho,
\eeq
then the system \eqref{LCEuler} can be equivalently written as the
isentropic flow for Chaplygin gas \cite{chap, serre} on varying cross-sectional area of the duct or with radially symmetry (see also equation (7.1.24) in \cite{Dafermos}), which is used for the modeling of dark energy,:
\beq
\label{LCEuler'}
\left\{
\begin{array}{l}
(c(x)\bar\rho)_t+(c(x)\bar\rho u)_x=0\\
(c(x)\bar\rho u)_t+\bigl(c(x)\bar\rho u^2- c(x)\bar\rho^{-1}\bigr)_x=- c'(x) \bar\rho^{-1}
\end{array}
\right.
\eeq
with pressure 
\[
p(\bar\rho)=-\bar\rho^{-1},
\]
where $\bar{\rho}$ is the density of gas at any point $(x,t)$ and $\rho$ is the density on a cross-section.
For a duct flow, $c(x)$ is the cross-sectional area which is uniformly positive and bounded.
We can also generally consider (\ref{LCEuler}) as a model 
for the isentropic Chaplygin gas in an inhomogeneous medium.

Similar as Theorem \ref{main-elast}, we construct blowups for the Cauchy problem of \eqref{LCEuler'} (or equivalently \eqref{LCEuler} by the relation \eqref{rho}).
For simplicity, we only consider a very special example of $c(x)$:
\beq
\label{c1}
 c(x)=\left\{
\begin{array}{rcl}
3-\e^\alpha -\frac{1}{2}\eta\e^\alpha,&x\in(-\infty,-1-\eta),\\
\psi_1(x),& x\in[-1-\eta,-1),\\
3+\e^\alpha x, & x\in[-1,1],\\
\psi_2(x),& x\in(1,1+\eta],\\
3+\e^\alpha +\frac{1}{2}\eta\e^\alpha,&x\in(1+\eta,\infty),\\
\end{array}\right.
\eeq
where $\alpha$ is any constant in $[0,1)$, $\psi_1(x)$ is an increasing function on $x\in[-1-\eta,-1)$ connecting $3-\e^\alpha -\frac{1}{2}\eta\e^\alpha$ and $3-\e^\alpha $ and $\psi_2(x)$ is an increasing function on $x\in[1,1+\eta]$ connecting $3+\e^\alpha $ and $3+\e^\alpha +\frac{1}{2}\eta\e^\alpha$. The positive constant $\eta<\e^3$. Furthermore $0< \e<1$ is a small given number which will be provided in the proof of the theorem.
By standard mollifier theory, we can find $\psi_1$ and $\psi_2$ such that $c(x)$ satisfies \eqref{c_con}. 

Now we list our first result on the singularity formation for the Cauchy problem of \eqref{LCEuler'}.

\begin{theorem}
\label{main}
For any $\alpha\in[0,1)$ and $c(x)$ given in \eqref{c1}, there exist uniformly bounded $C^1$ initial data $\rho_0(x)$ ($\rho_0\not\equiv1$ and has uniformly positive lower bound) and $u_0(x)$ (which will be given in the proof) and a finite time $T=O(\e^{-\alpha})>0$, such that the Cauchy problem of \eqref{LCEuler} with initial data $(\rho_0, u_0)$ has a unique classical $C^1$ solution $(\rho(x,t), u(x,t))$ on $(x,t)\in\mathbb{R}\times [0,T)$. Moreover, there exists a point 
$(x^*, T)$ with $|x^*|=O(\e^{\frac{1-\alpha}{2}})<O(1)$ at which the solution satisfies
\beq
\lim_{t\rightarrow{T^-}}u(x^*,t)=\infty \quad \lim_{t\rightarrow{T^-}}\rho(x^*,t)=0 \quad  \lim_{t\rightarrow{T^-}}\rho u(x^*,t)=B\quad \text{and}\quad
\lim_{t\rightarrow{T^-}}E(x^*,t)=\infty,
\eeq
where $B$ is a finite constant. 
The solutions $\rho(x,t)$ and $-u(x,t)$ have uniform upper bounds 
on $(x,t)\in\mathbb{R}\times [0,T)$.
\end{theorem}

\begin{remark}\label{remark1}
We have several remarks for Theorem \ref{main}.
\begin{itemize}
\item[(1).] In this example, the characteristic speeds of two families are uniformly away from each other,
or in another word, system \eqref{LCEuler'} is uniformly strictly hyperbolic (see Lemma \ref{lemmaR} and Lemma \ref{lemmaS}). 

\item[(2).] The result is also motivated by the pioneer work by Glassey, Hunter and Zheng \cite{ghz} on the singularities of variational wave equation \eqref{sn6} in one space dimension. Although equations \eqref{LCEuler'} and \eqref{sn6} are not equivalent when the initial $\rho_0\neq 1$,
initial data in Theorem \ref{main} and in the example in \cite{ghz} for \eqref{sn6} have a lot of similarities.


\item[(3).] This blowup can happen on a very slowly varying duct, which means $\|c'\|_{L^\infty}$ and the total variation of $c$ can be both arbitrarily small  in Theorem \ref{main}. When the variation of $c$ is larger ($c'(x)$  is larger)
around $x=0$, we show faster blowup. In fact, when $\alpha$ is decreasing, $c'(x)$ is increasing
around $x=0$, then the blowup time is shorter. When $\alpha=0$, the blowup time is at most $O(1)$.

\item[(4).] For any positive constant $K$, the total energy of the solution for any time until the blowup when $x\in[-K,K]$ is bounded by a constant depending on $K$, although at the time of blowup, the energy concentrates, i.e. the energy density is infinity, somewhere.

\end{itemize}
\end{remark}

\bigskip
When one looks for the radially symmetric solutions: $\bar \rho({\bf y},t)=\bar\rho(x,t)$, ${\bf u}({\bf y},t)={\bf y}u (x,t)$ with radius $x\geq0$ for
\beq\label{3deuler}
\left\{
\begin{array}{l}
\partial_t \bar\rho +\nabla_{\bf y}\cdot(\bar\rho {\bf u})=0\\
\partial_t (\bar\rho {\bf u}) +\nabla_{\bf y}\cdot(\bar\rho {\bf u}~{{\otimes}}~{\bf u}) +\nabla_{\bf y} p(\bar\rho)=0,
\end{array}
\right.
\eeq
with $p({\bf y},t)=\bar\rho^{-1}({\bf y},t)$, ${\bf u}({\bf y},t)=(u_1,u_2,u_3)$ and $({\bf y},t)\in\mathbb R^{m+1}\times \mathbb R^+$ with $m=1$ or $2$, the resulting system was in form of \eqref{LCEuler'} with 
$c(x)=x^m$ ($m=1$, cylindrical symmetric solution; $m=2$, spherically symmetric solution).
See \cite{Dafermos}. It will be shown in Section 1.3.1 that \eqref{3deuler} has strictly convex entropy.

Our next result concerns radially symmetric solutions for \eqref{3deuler}, that is the equation \eqref{LCEuler'} with $c(x)=x^m$, where $x$ denotes the radius.
Without of loss of generality, we only consider the solutions with initial data given on a special interval $x\in [1,3]$. 

\begin{theorem}
\label{main2} For $m=1,2$ and some given $C^1$ initial data $(\rho_0(x), u_0(x))$ depending only on radius $x\in[1,3]$ (which will be prescribed in the proof), there exists a time $T= O(1)>0$ such that the Cauchy problem of \eqref{LCEuler'} with $c(x)=x^m$ has a unique classical $C^1$ solution in $\Omega_{symm}$, where  $\Omega_{symm}$ is the domain of dependence of the initial interval  $x\in[1,3]$ for any time $t$ in $(0,T)$. Moreover, there exists a point 
$(x^*, T)$ with $|x^*|=O(\e^{\frac{1}{2}})<O(1)$ at which the solution satisfies
\beq
\lim_{t\rightarrow{T^-}}u(x^*,t)=\infty \quad \lim_{t\rightarrow{T^-}}\rho(x^*,t)=0 \quad  
\lim_{t\rightarrow{T^-}}\rho u(x^*,t)=B\quad \text{and}\quad
\lim_{t\rightarrow{T^-}}E(x^*,t)=\infty,
\eeq
where $B$ is a finite constant. 
The solutions $\rho(x,t)$ and $-u(x,t)$ have uniform upper bounds 
on $(x,t)\in\mathbb{R}\times [0,T)$.
\end{theorem}

\begin{remark}
Theorem \ref{main2} provides an example with finite time vacuum formation and $L^\infty$ blowup
for the radially symmetric solutions with radius varying in a finite closed interval away from zero. 
This example satisfies all properties as the example in Theorem \ref{main} with $\alpha=0$.

\end{remark}



%
\subsection{Convex entropy and vacuum}
Smooth solutions of the system \eqref{3deuler} satisfy energy conservation law
\[
{\mathcal E}_t+\nabla_{\bf y}\cdot {\bf Q}=0
\]
with entropy
\[
\mathcal E=\frac{1}{2}\bar\rho |{\bf u}|^2+ \frac{1}{2}\bar\rho^{-1}
\]
and entropy flux
\[
{\bf Q}=(\frac{1}{2}\bar\rho |{\bf u}|^2-\frac{1}{2}\bar\rho^{-1}){\bf u}.
\]
The entropy $\mathcal{E}$ is a strictly convex function on conservative variables $(\bar\rho,{\bf m})$,
where ${\bf m}=\bar\rho{\bf u}$ is the momentum,
i.e. 
$\mathcal E=\frac{1}{2}\frac{|{\bf m}|^2}{\bar\rho}+ \frac{1}{2}\bar\rho^{-1}$
satisfies that $D^2 \mathcal E$ is a positively defined matrix.

Concerning one dimensional case, the smooth solutions of equation {\eqref{LCEuler}} satisfy energy conservation law
$$
E_t+q_x=0,
$$
with entropy
$$
E=\frac{1}{2}\rho u^2+\frac{1}{2}c^2\rho^{-1}
$$
and entropy flux
$$
q=\frac{1}{2}u^3\rho-\frac{1}{2}c^2\rho^{-1}u.
$$
By direct calculation, we obtain
$$
E=\frac{1}{2}
\left(\begin{array}{lcr}
u&\rho
\end{array}
\right)
\left(\begin{array}{lcr}
\rho&0\\
0&c^2\rho^{-3}
\end{array}
\right)
\left(\begin{array}{lcr}
u\\
\rho
\end{array}
\right).
$$ 
This implies that the entropy $E$ is strictly convex.

The $L^\infty$ blowup in this paper is totally different from the previous $L^\infty$ blowup results found first by Jenssen in his groundbreaking work \cite{je}, and then by several other authors \cite{bj,jy,yng,ys} by considering the shock interactions. To see the difference, more intuitively, one could still essentially consider that
the $L^\infty$ blowup constructed in our examples on $\frac{1}{\rho}$ are coming from the blowup on the gradient variable $x_X$ through transformation between different coordinates. A key point worth mentioning is that after transformation from Lagrangian coordinates to Eulerian coordinates, one gets a linearly degenerate system, which {\em rarely} happen, on which $C^1$ solution exists before the $L^\infty$ blowup. 

Furthermore, except the $L^\infty$ blowup, the more generic singularity: gradient blowup has been studied for systems of conservation laws. The gradient blowup in systems of conservation laws has been widely accepted as the most generic type of singularity related to the shock formation in \cite{lax,john,G3, G5, G6, G8}. It is much harder to find the $L^\infty$ blowup for systems of conservation laws.


%
Finally, we give a remark on the  $L^\infty$ blowup on $u$ and vacuum formation. The isentropic hyperbolic systems
with $p={\bar \rho}^\gamma$ when the adiabatic constant $\gamma>0$ and  $p=-{\bar \rho}^\gamma$ when $\gamma<0$ are given by
\beq
\label{LCEuler-1}
\left\{
\begin{array}{l}
(c(x)\bar\rho)_t+(c(x)\bar\rho u)_x=0\\
(c(x)\bar\rho u)_t+\bigl(c(x)\bar\rho u^2+c(x)p\bigr)_x= c'(x) p.
\end{array}
\right.
\eeq
When $0\leq\gamma< 1$, the entropy is not strictly convex so these cases are not the interesting cases for us.
When $\gamma\geq 1$ or $-1<\gamma<0$ or $\gamma<-1$, the system \eqref{LCEuler-1} is genuinely nonlinear when $c(x)$
is a constant, hence we expect shock formation and
tend to believe
that the shock could prevent the $L^\infty$ blowup.
For example, when $1<\gamma\leq\frac{5}{3}$, which is corresponding to gas dynamics,
the $L^\infty$ existence has already been proven in \cite{GG} for the duct flow and exterior radially symmetric flow, hence $L^\infty$ blowup on $u$ 
can not happen (see also \cite{lm} for the gas dynamics with $1<\gamma<\infty$). 

Furthermore, the nonisentropic compressible Euler equations with polytropic ideal gas are
\beq\label{Euler}
\left\{\begin{array}{l}
  \rho_t +(\rho\,u)_{x} = 0\\
  (\rho\,u)_t + (\rho\,u^2 + p)_{x} = 0 \\ 
  (\frac{1}{2}\rho\,u^2 + \rho\,e)_t + (\frac{1}{2}\,\rho\,u^3 + u\,p)_{x} = 0\,,
\end{array}
\right.
\eeq
with equation of state
\[
 e=c_v T={\frac{p\,\tau}{\gamma-1}} \quad\text{and}\quad p\,\tau=\mathcal{R}\,T,
\]
so that pressure
\beq
  p=K\exp(\mathcal{S}/c_v)\,\tau^{-\gamma},\label{poftau}
\eeq 
where $\rho$ is density, $\tau=1/\rho$, $u$ is velocity, $e$ is specific 
internal energy, $\mathcal{S}$ is the entropy, $T$ is the temperature, $\mathcal{R}$, $K$, $c_v$ are
positive constants, and  $\gamma>1$ is the adiabatic gas constant.
In \cite{G8}, the first author, R. Young and Q. Zhang have found 
uniform time-independent $L^\infty$ bounds for $\rho$ and $|u|$.  

Whether the solution for \eqref{LCEuler-1} or \eqref{Euler} with $\gamma>1$ has a finite time 
vacuum or not is still a major open problem for gas dynamics. 
If one only
considers the smooth solutions for \eqref{LCEuler-1} with $\gamma>1$ and constant
$c$, one may conjecture, by a strong evidence from \cite{lin}, that there will be no vacuum 
in finite time if there is no vacuum initially or instantaneously.









\bigskip
The rest of paper is organized as follows.  
In Section 2, we set up the Riemann coordinates and several lemmas for smooth solutions.
In Section 3, we prove the existence of $C^1$ solution when 
$\|u\|_{L^{\infty}}+\|\rho\|_{L^{\infty}}<+\infty$ and $|\rho|$ is away from zero.
In Section 4, we prove Theorem \ref{main-elast} for inhomogeneous elastic flow.
In Section 5, we prove the Theorem \ref{main} for isentropic Chaplygin gas and the Theorem \ref{main2} for the radially symmetric case.

\section{Riemann coordinates}
\setcounter{equation}{0}
\setcounter{theorem}{0}

For smooth solutions, the system \eqref{LCEuler} can be written as
\begin{equation}
\label{vac-eq2}
\begin{split}
\rho_t+u\rho_x +\rho u_x&=0\\
u_t+uu_x+c^2(x)\rho^{-3}\rho_x&=c(x)c'(x)\rho^{-2}.
\end{split}
\end{equation}
Hence
\beq\left(
\begin{array}{c}
\rho\\
u
 \end{array}
\right)_t
+
A
\left(
\begin{array}{c}
\rho\\
u
 \end{array}
\right)_x
=
\left(
\begin{array}{c}
0\\
c(x)c'(x)\rho^{-2}
 \end{array}
\right)
\eeq
where
\begin{equation}
\label{vac-eq3}
A=\left(
\begin{array}{cc}
u& \rho\\
c^2\rho^{-3}& u
\end{array}
\right).
\end{equation}
Direct calculation shows that the eigenvalues of $A$ are 
\begin{equation}
\label{vac-eq4}
S=u+c(x)\rho^{-1},\quad R=u-c(x)\rho^{-1},
\end{equation}
and the corresponding right eigenvectors are 
\begin{equation}
\label{vac-eq5}
v_1=(1,\ c(x)\rho^{-2})^T,\quad v_{2}=(1,\ -c(x)\rho^{-2})^T.
\end{equation}
According to Lax \cite{lax0}, the two characteristic families for system \eqref{vac-eq2} when $c$ is a constant are both {\em linearly degenerate}.
For any $(\bar{x},\bar{t})$ with $\bar{t}>0$, the plus and minus characteristics 
$x_{\pm}(t)$ through $(\bar{x}, \bar{t})$ are defined by 
\begin{equation}
\label{vac-eq17}
\frac{dx_{+}(t,\bar{x}, \bar{t})}{dt}=S(x_{+},t)\qquad \text{and}\qquad 
\frac{dx_{-}(t, \bar{x}, \bar{t})}{dt}=R(x_{-},t)\,. 
\end{equation}
For simplicity, we use $x_+(t)$ or $t_+(x)$ and $x_-(t)$ or $t_-(x)$
to denote the plus and minus characteristics respectively.
For smooth solutions of equation \eqref{vac-eq2}, we obtain the following equation of $c(x)/\rho$, which will be used several times in the rest of paper.
\begin{equation}
\label{vac-eq6}
c(x)(\rho^{-1})_t+u(c(x)\rho^{-1})_x-c(x)\rho^{-1}u_x=c'(x)u\rho^{-1}.
\end{equation}

\begin{lemma}
\label{lemma2.1} 
For smooth solutions of \eqref{LCEuler}, $S$ and $R$ satisfy 
\begin{equation}
\label{vac-eq12}
\begin{split}
S_t+RS_x=c'(x)u\rho^{-1}=\frac{c'(x)}{4c(x)}\left(S^2-R^2\right),\\
R_t+SR_x=-c'(x)u\rho^{-1}=\frac{c'(x)}{4c(x)}\left(R^2-S^2\right).
\end{split}
\end{equation}
\end{lemma}
\pf
By equation \eqref{vac-eq2} and \eqref{vac-eq6}, we have 
\begin{equation}
\label{vac-eq7}
\begin{split}
&S_t+RS_x\\
=&(u+c(x)\rho^{-1})_t+(u-c(x)\rho^{-1})(u+c(x)\rho^{-1})_x\\
=&u_t+c(x)(\rho^{-1})_t+uu_x-c(x)\rho^{-1}u_x+u(c(x)\rho^{-1})_x-c(x)\rho^{-1}(c(x)\rho^{-1})_x\\
=&u_t+uu_x-c(x)\rho^{-1}(c(x)\rho^{-1})_x+c(x)(\rho^{-1})_t-c(x)\rho^{-1}u_x+u(c(x)\rho^{-1})_x\\
=&c'(x)u\rho^{-1}.
\end{split}
\end{equation}
And
\begin{equation}
\label{vac-eq8}
\begin{split}
&R_t+SR_x\\
=&(u-c(x)\rho^{-1})_t+(u+c(x)\rho^{-1})(u-c(x)\rho^{-1})_x\\
=&u_t-c(x)(\rho^{-1})_t+uu_x+c(x)\rho^{-1}u_x-u(c(x)\rho^{-1})_x-c(x)\rho^{-1}(c(x)\rho^{-1})_x\\
=&u_t+uu_x-c(x)\rho^{-1}(c(x)\rho^{-1})_x-\left(c(x)(\rho^{-1})_t-c(x)\rho^{-1}u_x+u(c(x)\rho^{-1})_x\right)\\
=&-c'(x)u\rho^{-1}.
\end{split}
\end{equation}
\endpf

\begin{remark}
\label{remark2.2}
When $c(x)$ is a constant function, $R$ and $S$ are exactly two Riemann invariants
along plus and minus characteristics with characteristic speeds $S$ and $R$ respectively by Lemma \ref{lemma2.1}, 
\beq\label{cx0}
\left\{\begin{array}{rcl}
S_t+RS_x&=&0,\vspace{.2cm}\\
R_t+SR_x&=&0.
\end{array}
\right.
\eeq
By the Theorem 2.3 in \cite{li1}, system (\ref{cx0}) admits a global-in-time unique $C^1$ solution
if the initial data $S_0$ and $R_0$ have bounded $C^1$ norm.

\end{remark}

\begin{lemma}[Energy conservation law] 
For smooth solutions of \eqref{LCEuler}, the energy density $E=\frac{1}{4}\rho(S^2+R^2)$ satisfies
\begin{equation}\label{vac-eq9}
\left(\rho(S^2+R^2)\right)_t+\left(\rho(S^2R+R^2S)\right)_x=0.
\end{equation}
\end{lemma}
\pf
By\eqref{vac-eq4}, we have
\begin{equation}\label{vac-eq10}
\rho S=\rho u+c(x)\quad\text{and}\quad \rho R=\rho u-c(x).
\end{equation}
Thus, multiplying the first equation of \eqref{vac-eq12} by $2\rho S$ 
and the second equation of \eqref{vac-eq12} by $2\rho R$, adding them up, and using \eqref{vac-eq10} and the conservation of mass in \eqref{LCEuler}, 
we obtain
\begin{equation}
\label{vac-eq11}
\begin{split}
&(\rho S^2+\rho R^2)_t+\left(\rho(S^2R+R^2S)\right)_x\\
=&(S^2+R^2)\rho_t+S^2(\rho R)_x+R^2(\rho S)_x+2c'u(S-R)\\
=&(S^2+R^2)\rho_t+S^2(\rho u-c)_x+R^2(\rho u+c)_x+c'(S+R)(S-R)\\
=&S^2(\rho_t+(\rho u)_x)+R^2(\rho_t+(\rho u)_x)+c'(R^2-S^2)+c'(S^2-R^2)\\
=&0.
\end{split}
\end{equation}
\endpf

\begin{figure}[htb]
\centering
\includegraphics[width=7cm,height=4cm]{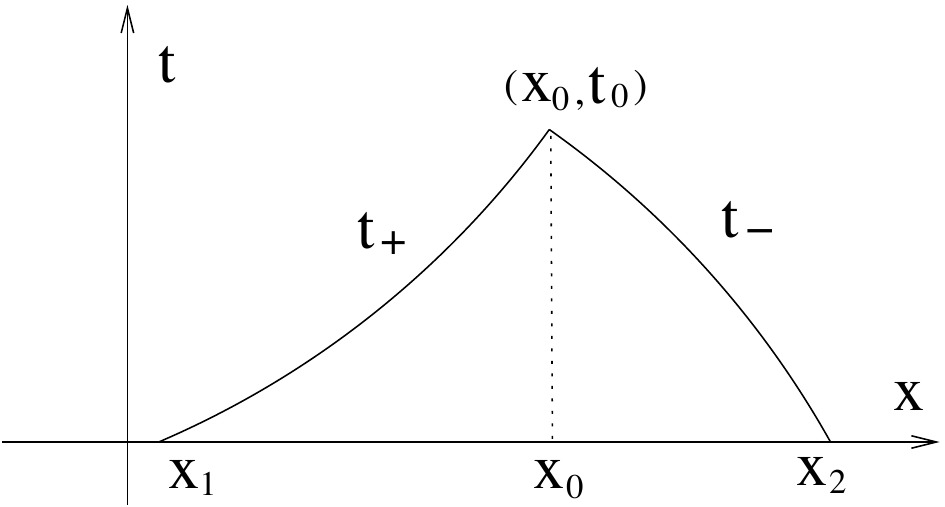}
\caption{A characteristic triangle $D$.}
\label{f0}
\end{figure}

Finally, we give a key estimate for the proof of our main theorems.
For any $(x_0,t_0)\in\mathbb R\times(0,+\infty)$, let $\gamma_{+}$ and $\gamma_{-}$ be plus and minus characteristics through $(x_0,t_0)$ and intersect $x-$axis at $x_1$ and $x_2$ respectively (see Figure \ref{f0}).
Integrating \eqref{vac-eq9} over a characteristic triangle $D$ enclosed by $\gamma_{+}$, $\gamma_{-}$ and $[x_1,x_2]$ as in Figure \ref{f0}, 
we obtain an energy identity indicating the finite propagation of the waves.
\begin{lemma}[Finite propagation]
\label{vac-lemma1}
For smooth solutions of \eqref{LCEuler} inside a characteristic triangle $D$ in Figure \ref{f0},
\begin{equation}
\label{vac-eq18}
2\int_{x_1}^{x_0}cS(x,\gamma_{+}(x))\,dx-2\int^{x_2}_{x_0}cR(x,\gamma_{-}(x))\,dx
=\int_{x_1}^{ x_2}\rho(S^2+R^2)(x,0)\,dx.
\end{equation}
\end{lemma}

\pf Integrating \eqref{vac-eq9} over the region $D$ and using the Green's theorem, we obtain
\begin{equation}
\label{vac-eq19}
\begin{split}
0=&\iint_{D}\left(\rho(S^2+R^2)\right)_t+\left(\rho(S^2R+R^2S)\right)_x\,dxdt\\
=&\int_{\partial D}\rho(S^2R+R^2S)\,dt-\rho(S^2+R^2)\,dx.
\end{split}
\end{equation}
On $\gamma_{+}$ which is the left boundary of $D$, by \eqref{vac-eq4} and \eqref{vac-eq17} we have
\begin{equation}
\label{vac-eq20}
\begin{split}
&\int_{\gamma_{+}}\rho(S^2R+R^2S)\,dt-\rho(S^2+R^2)\,dx\\
=&\int_{x_1}^{x_0}\rho(SR+R^2)-\rho(S^2+R^2)\,dx\\
=&\int_{x_1}^{x_0}\rho S(R-S)\,dx\\
=&-2\int_{x_1}^{x_0}cS(x,\gamma_{+}(x))\,dx.
\end{split}
\end{equation}
Similarly, on $\gamma_{-}$ which is the right boundary of $D$, we have
\begin{equation}
\label{vac-eq21}
\begin{split}
&\int_{\gamma_{-}}\rho(S^2R+R^2S)\,dt-\rho(S^2+R^2)\,dx\\
=&\int^{x_2}_{x_0}\rho(S^2+RS)-\rho(S^2+R^2)\,dx\\
=&\int^{x_2}_{x_0}\rho R(S-R)\,dx\\
=&2\int^{x_2}_{x_0}cR(x,\gamma_{-}(x))\,dx.
\end{split}
\end{equation}
Putting \eqref{vac-eq20} and \eqref{vac-eq21} into \eqref{vac-eq19}, we obtain \eqref{vac-eq18}.
\endpf

%

\section{Existence of $C^1$ solutions before $L^{\infty}$ blowup}
\setcounter{equation}{0}
\setcounter{theorem}{0}
For smooth solutions of \eqref{LCEuler}, let 
$$v=\rho^{-1}S_x,\qquad w=\rho^{-1}R_x$$
be two gradient variables.
Then we obtain the following lemma.
\begin{lemma} 
\label{lemma3.1}
For smooth solutions of \eqref{LCEuler}, $v$ and $w$ satisfy
\begin{equation}
\label{vac-eq13}
\begin{split}
v_t+Rv_x=&\frac{c'}{2c}\left[(2c\rho^{-1}+S)v
-Rw\right]
+\frac{c''c-c'^2}{4c^2}\rho^{-1}(S^2-R^2)
\\
w_t+Sw_x=&\frac{c'}{2c}\left[(2c\rho^{-1}+R)w
-Sv\right]
+\frac{c''c-c'^2}{4c^2}\rho^{-1}(R^2-S^2)\,.
\end{split}
\end{equation}
\end{lemma}

\pf 
By \eqref{vac-eq12} and \eqref{vac-eq4}, we have
\begin{equation}
\label{vac-eq14}
\begin{split}
&v_t+Rv_x\\
=&S_x\left(\rho^{-1}\right)_t+S_xR\left(\rho^{-1}\right)_x
+\rho^{-1}\left(S_x\right)_t+\rho^{-1}R\left(S_x\right)_x\\
=&S_x\left(\left(\rho^{-1}\right)_t+u\left(\rho^{-1}\right)_x-c\rho^{-1}\left(\rho^{-1}\right)_x\right)
+\rho^{-1}\left(S_t+RS_x\right)_x-\rho^{-1}R_xS_x\\
=&S_x\left(u_x\rho^{-1}-c\rho^{-1}\left(\rho^{-1}\right)_x\right)
+\rho^{-1}\left(S_t+RS_x\right)_x-\rho^{-1}R_xS_x\\
=&S_x\rho^{-1}\left(u_x-\left(c\rho^{-1}\right)_x+c'\rho^{-1}\right)
+\rho^{-1}\left(S_t+RS_x\right)_x-\rho^{-1}R_xS_x\\
=&\rho^{-1}R_xS_x+c'\rho^{-2}S_x
+\rho^{-1}\left(S_t+RS_x\right)_x-\rho^{-1}R_xS_x\\
=&c'\rho^{-2}S_x
+\rho^{-1}\left(S_t+RS_x\right)_x\\
=&c'\rho^{-2}S_x
+\rho^{-1}\left(\frac{c'(x)}{4c(x)}\left(S^2-R^2\right)\right)_x,
\end{split}
\end{equation}
which implies the first equation of \eqref{vac-eq13}.
Similarly
\begin{equation}
\label{vac-eq15}
\begin{split}
&\left(\rho^{-1}R_x\right)_t+S\left(\rho^{-1}R_x\right)_x\\
=&R_x\left(\rho^{-1}\right)_t+R_xS\left(\rho^{-1}\right)_x
+\rho^{-1}\left(R_x\right)_t+\rho^{-1}S\left(R_x\right)_x\\
=&R_x\left(\left(\rho^{-1}\right)_t+u\left(\rho^{-1}\right)_x-c\rho^{-1}\left(\rho^{-1}\right)_x\right)
+\rho^{-1}\left(R_t+SR_x\right)_x-\rho^{-1}R_xS_x\\
=&R_x\left(u_x\rho^{-1}-c\rho^{-1}\left(\rho^{-1}\right)_x\right)
+\rho^{-1}\left(R_t+SR_x\right)_x-\rho^{-1}R_xS_x\\
=&R_x\rho^{-1}\left(u_x-\left(c\rho^{-1}\right)_x+c'\rho^{-1}\right)
+\rho^{-1}\left(R_t+SR_x\right)_x-\rho^{-1}R_xS_x\\
=&\rho^{-1}R_xS_x+c'\rho^{-2}R_x
+\rho^{-1}\left(R_t+SR_x\right)_x-\rho^{-1}R_xS_x\\
=&c'\rho^{-2}R_x
+\rho^{-1}\left(R_t+SR_x\right)_x\\
=&c'\rho^{-2}R_x
+\rho^{-1}\left(\frac{c'(x)}{4c(x)}\left(R^2-S^2\right)\right)_x,
\end{split}
\end{equation}
which implies the second equation of \eqref{vac-eq13}.
{\endpf}

\begin{remark}\label{remark3.2}
System \eqref{vac-eq13} indicates that the rates of change of $v$ along the minus characteristic and $w$ along the plus characteristic are both linear.
\end{remark}


Before we prove the $C^1$ existence result when the solution has $L^\infty$ bounds, 
we first state an a priori condition.

\begin{itemize}
\item[{\bf (A)}] Suppose the initial data $u_0$ and $\rho_0$ are $C^1$ functions and 
uniformly bounded ($\rho_0$ is also uniformly bounded away from zero). Then for
any $C^1$ solution $(\rho(x,t),~ u(x,t))$ of Cauchy problem of equation 
\eqref{LCEuler} with $(x,t)\in \mathbb{R}\times[0,T_*]$, for $0<T_*<T$,
there exists a positive constant $L_*$, 
only depending on $T_*$ and initial data, such that
\beq\label{L_def}
\|\rho(x,t)\|_{L^{\infty}}+\|\rho^{-1}(x,t)\|_{L^{\infty}}+\|u(x,t)\|_{L^{\infty}}=L_*<\infty.
\eeq 

\end{itemize}
Under the condition (A), by Lemma \ref{lemma3.1}, observation in Remark \ref{remark3.2}
and $C^{\infty}$ functions are dense in $C^1$, it is easy to obtain the following a priori estimates.
\begin{lemma}
\label{lemma3.3}
Assume that the condition (A) is satisfied, $u_0, \rho_0\in C^1$ and uniformly bounded ($\rho_0$ is also uniformly bounded away from zero). Then any $C^1$ solutions $\rho(x,t)$, $u(x,t)$ for system \eqref{LCEuler} with $0\leq t\leq T_*$ satisfy 
\begin{equation}
\label{est-1st-deriv}
\sup_{(x,t)\in \mathbb{R}\times[0,T^*]}\{|S_t|, |S_x|, |R_t|, |R_x|\}=M_*,
\end{equation}
for some positive constant $M_*$ only depending on initial values, $L_*$ and $T_*$.
\end{lemma}

\begin{remark}
Lemma \ref{lemma3.3} implies that the $C^1$ blowups never occur before $L^{\infty}$ blowup for \eqref{LCEuler}. 
\end{remark}

By the a priori estimate in Lemma \ref{lemma3.3}, one can prove the existence of $C^1$ solution on $[0,T)$ under the condition (A). 

\begin{theorem}
\label{ex-sm-sl}
If the condition (A) is satisfied for any $0<T_*<T$,
the initial value problem of \eqref{LCEuler} with uniformly bounded initial data 
$u_0, \rho_0\in C^1$ ($\rho_0$ is also uniformly bounded away from zero) has a unique 
$C^1$ solution $(\rho(x,t),~ u(x,t))$ for $(x,t)\in \mathbb{R}\times[0,T)$.  
\end{theorem}

First we observe that (A) implies the uniformly strict hyperbolicity of \eqref{LCEuler}. The local existence of the $C^1$ solution now can be obtained by standard argument in \cite{Liyu}.  
Under the condition $(A)$, the local solution can be extended to $(x,t)\in \mathbb{R}\times[0,T)$ 
by Theorem 2.4 in \cite{li1} and Remark 2.20 in \cite{li1}. 
To make this paper self-contained, we sketch the proof here. 

\pf
We first fix our consideration on $(x,t)\in \mathbb{R}\times[0,T_*]$ for some $T_*$.

By the local-in-time existence results for the quasi-linear first order hyperbolic systems 
in \cite{Liyu},
for any strong determinate domain $\tilde\Omega_k$ corresponding to initial interval $x\in[-k,k]$,
there exists some time $T_0=T_0(k, L_*, M_*)$  
such that \eqref{vac-eq13} exists a unique $C^1$ 
solutions on $\tilde\Omega_k$ with $t\in[0,T_0]$. Here a domain 
\[
\tilde\Omega_{[a,b]}\equiv\tilde\Omega_{[a,b]}(\delta_0)=\left\{(x,t)\big| 0\leq t\leq \delta_0,\ x_1(t)\leq x\leq x_2(t)  \right\}
\] 
is called a strong determinate domain of initial interval $[a,b]$ if
\begin{itemize}
 \item[i.] $x_1(t)$ and $x_2(t)$ are $C^1$ functions for  $ 0\leq t\leq \delta_0$.
 \item[ii.] $x_1(0)=a$ and $x_2(0)=b$.
 \item[iii.] For any $C^1$ solution in $\tilde\Omega(\delta_0)$, $x'_1(t)\geq M_S$ and $x'_2(t)\leq M_R$,
\end{itemize}
where $M_S$ is the upper bound of plus characteristic $S$ and 
$M_R$ is the lower bound of minus characteristic $R$ on $\tilde\Omega(\delta_0)$. In our problem,
$M_S$ and $M_R$ are only depending on $L_*$.
 Since $T_0$ is a constant,  when $k$ is large enough, one can prove the existence of $C^1$ solution on $\tilde\Omega_k$ with $t\in[0,T_*]$,
by using the local existence proof finite many times.

Next, for any point $(x,t)\in \mathbb{R}\times[0,T_*]$, we could find  $\tilde\Omega_k$ including this point
with sufficiently large $k$,  because $M_S$ and $M_R$ are only depending on $L_*$. 
Hence, we already proved the global existence on $(x,t)\in \mathbb{R}\times[0,T_*]$.

Furthermore, since $T_*$ is any time before $T$, so we already proved the $C^1$ existence on $(x,t)\in \mathbb{R}\times[0,T)$,
where the uniqueness of the local existence protects that we have one unique solution.
\endpf
\section{Vacuum for the inhomogeneous elastic flow: proof of Theorem \ref{main-elast}}
\setcounter{equation}{0}

We define an increasing $C^2$ smooth function
\beq
\label{elas1_c0}
 c(x)=\left\{
\begin{array}{rcl}
d_1,&x\in(-\infty,-\e^{-8}-\e^5),\\
\Psi_1(x),& x\in[-\e^{-8}-\e^5,-\e^{-8}),\\
\frac{\e}{1-\e x}, & x\in[-\e^{-8},1],\\
\Psi_2(x),& x\in(1,1+\e^{5}],\\
d_2,&x\in(1+\e^{5},\infty),\\
\end{array}\right.
\eeq
where $d_1$, $d_2$ are two positive constants and $\Psi_1$, $\Psi_2$ are increasing $C^2$ smooth functions
such that 
\beq\label{elas_ineq_d}
\frac{\e}{1-\e}\leq d_2\leq \frac{\e}{1-\e}+\e^5,
\eeq
and
\beq\label{elas_ineq_d2}
\frac{\e}{1+\e(\e^{-8}+\e^5)}(1-\e^5)\leq d_1\leq \frac{\e}{1+\e(\e^{-8}+\e^5)}.
\eeq
The function  $f(x)=\frac{\e}{1-\e x}$ on $x\in[-\e^{-8},1]$ satisfies
\[
\frac{f'(x)}{f^2(x)}=1,
\]
hence we can find $d_1$, $d_2$, $\Psi_1(x)$ and $\Psi_2(x)$ such that
\beq\label{elas_ineq_c}
0\leq\frac{c'(x)}{c^2(x)}\leq 1+ \e^8
\eeq
for any $x$. It is easy to see that there exists a function  
$c(x)$ such that \eqref{c_con}, \eqref{elas1_c0}$\sim$\eqref{elas_ineq_c}
are all satisfied. 
The positive constant $\e\ll1$ will be given in the proof of the theorem.

\begin{remark} \label{elas_remarkRho0}
The construction of the initial data $u_0$ and $c(x)$ is motivated by the seminal work \cite{ghz}. However, our restrictions are all on the function $c(x)$ since
$\rho_0\equiv 1$, $u_0=c(x)$ and $S$ is uniformly larger than $R$. It makes
our construction much more involved than the one in \cite{ghz}. 
\end{remark}

Throughout this paper, we use $K_i$ and $M_i$ to denote positive constants
independent of $\e$.
To prove Theorem \ref{main-elast}, we show that there
exists some time 
 \beq\label{elas_T_def}
T=M_0 \e^{-2}
 \eeq
such that the a priori condition (A) is satisfied for any $t\in[0,T^*]$ with  $T^*<T$ which indicates that $C^1$ solution exists for any $t\in[0,T^*]$. We also show that at $T=M_0 \e^{-2}$, $S$ blows up at somewhere while $R$ is uniformly bounded.

\begin{lemma}\label{elas_lemma0}
For any $C^1$ solutions to \eqref{LCEuler} with $c(x)$ given in $\eqref{elas1_c0}$ and prescribed initial data $\rho_0\equiv1$, $u_0=c(x)$, we have
 \[
 R\leq0.
 \]
\end{lemma}
\pf
By  Lemma \ref{lemma2.1}, for any $x\in(-\infty,\infty)$
\beq\label{elas_Rin0}
R_t+S R_x \leq M_1(\e, d_1) R^2,
\eeq
where $M_1(\e,d_1)$ is a positive constant depending on $\e$ and $d_1$, and left hand side is 
the derivative along a minus characteristic.
By \eqref{elas_Rin0} and ODE comparison theorem, 
\[R\leq0\]
 for any $C^1$ smooth solution.
\endpf
\vspace{.2cm}

In Lemmas \ref{elas_lemmaT}$\sim$\ref{elas_lemmaS}, we restrict our consideration on the $C^1$ solutions with $S(x,t)>0$ for any $(x,t)\in\mathbb R\times [0,T^*]$ with $T^*<T$ and $T$ defined in \eqref{elas_T_def}. For these solutions, the plus and minus characteristics go in forward and backward directions respectively. 
We will show that the a priori condition on $S$ is satisfied for any $C^1$ solutions with prescribed initial data in Lemma \ref{elas_Sgeq0}.

\begin{figure}[htp]
\begin{center}
\includegraphics[width=12.5cm,height=5.8cm]{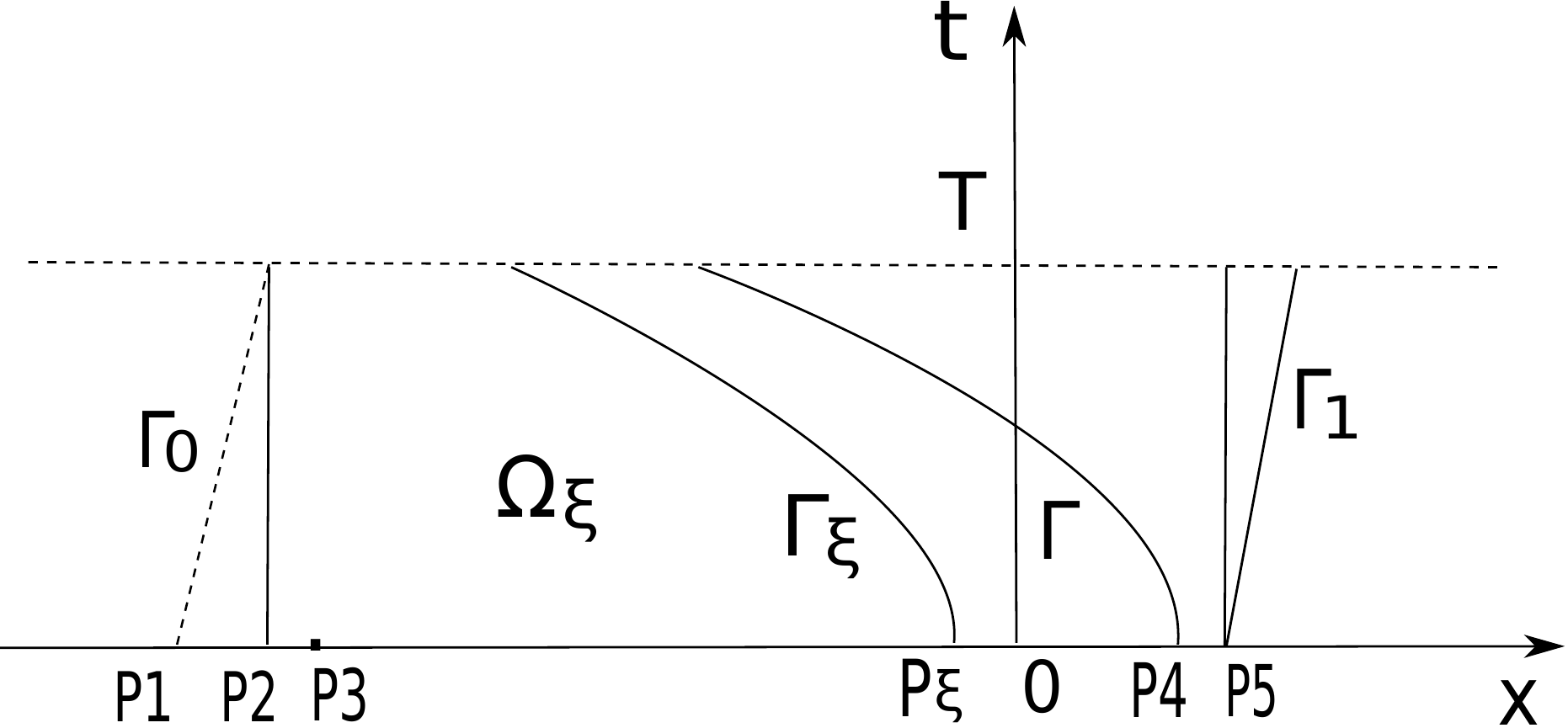}
\end{center}
{\caption{Proof of Theorem \ref{main-elast}. }}
\label{euler2}
\end{figure}

We use Figure \ref{euler2} for the proof.
Note in Figure \ref{euler2}, $P_1$ is the initial point of the forward
characteristic $\Gamma_0$ intersecting with the vertical line passing $P_2$ at time $T$,
$P_2=(-\e^{-8}-\e^5,0)$, $P_3=(-\e^{-8},0)$, $P_4=(1,0)$ and $P_5=(1+\e^{5},0)$.
The backward characteristic $\Gamma_\xi$ starts from $P_\xi=(\xi,0)$
and ends at $t=T$ with $x_{P_2}<\xi\leq1$. The backward characteristic $\Gamma$ starts from $P_4$. The forward characteristic $\Gamma_1$
starts from $P_5$.
Furthermore, we use $\Omega_{\xi}$ to denote the domain of dependence in the left of $\Gamma_\xi$ when $t<T$.

\begin{lemma}\label{elas_lemmaT}
Consider any $C^1$ solutions for \eqref{LCEuler} with $c(x)$ given in $\eqref{elas1_c0}$ and prescribed initial data $\rho_0\equiv1$, $u_0=c(x)$. Assume $S(x,t)>0$ for any $(x,t)\in\mathbb R\times [0,T)$. 
Then the initial energy $E_\xi=\frac{1}{4}\int_{P_1 P_\xi}\rho(S^2+R^2)(x,0)\,dx\ $ in the interval $P_1 P_\xi$ with $x_{P_1}<\xi\leq1$ satisfies
\beq\label{elas_e_xi}
0<(1+\e^8)E_\xi\leq \frac{\e}{1-\e\xi}\,.
\eeq
In the region $\Omega_{\xi}$,
\beq\label{elas_lemmaT_R}
0\geq R(x,t)\geq- \frac{\e}{2(1-\e\xi)}=:-a_\xi\,.
\eeq
Also we know that the characteristic $\Gamma$ will not interact with $\Gamma_0$
before $T$.
Furthermore,
\beq\label{elas_S_0}
S_0(\xi)=\frac{2\e}{1-\e\xi}(1-\e^5)>a_\xi\,,
\eeq
when $x_{P_2}<\xi\leq1$.
\end{lemma}
%
\pf
By Figure \ref{euler2} and the initial data, the energy in the initial interval $P_3 P_\xi$ is
\[
E_{P_3 P_\xi}=\int_{-\e^{-8}}^{\xi} \,\bigl(\frac{\e}{1-\e x}\bigr)^2~dx=\frac{\e}{1-\e \xi}-\frac{\e}{1-\e^{-7}}.
\]
Since the length of $P_2P_3$ is $\e^{5}$ and the forward characteristic $\Gamma_0$ has speed less than $O(\e^{7})$ which are both very small, so it is easy to see that energy $E_{P_1P_3}$ is omittable hence
(\ref{elas_e_xi}) is correct when $\e$ is small enough, where we also use that $T=O(1)\e^{-2}$.

For any $(x,t)$ in the left the forward characteristic $\Gamma_0$, 
$R(x,t)=0$ since $c(x)$ is a constant when $x<-\e^{-8}-\e^5$ and \eqref{cx0}, hence \eqref{elas_lemmaT_R} is correct in this region.

By  \eqref{vac-eq12}, we have
\beq\label{elas_R_ineq}
R_t+S R_x \geq - \frac{c'}{4c^2}S^2.
\eeq
Integrating it along any forward characteristic in $\Omega_{\xi}$ starting from $P_1 P_4$, by \eqref{vac-eq18} and $\frac{c'}{c^2}\leq 1+\e^{8}$, for any 
$(x,t)$ on this forward characteristic  in $\Omega_{\xi}$ we have
\beq
R(x,t)\geq - \int_{\gamma_+\cap\Omega} \frac{c'}{8c^2}\,2cS^2(t,x_{+}(t))dt\geq-\frac{1}{2}(1+\e^{8})E_\xi\geq- \frac{\e}{2(1-\e\xi)}=-a_\xi.\label{elas_R_est}
\eeq
Hence, \eqref{elas_lemmaT_R} is always correct in $\Omega_\xi$.

For any $(x,t)$ on $\Gamma$ with wave speed $R$, we have
\[
\frac{1-x}{t}\leq  \frac{\e}{2(1-\e)}=O(\e).
\]
Since,  $|P_1P_4|>\e^{-8}$, easy calculation shows that the characteristic $\Gamma$ will not interact with $\Gamma_0$
before $T=O(\e^{-2})$.

Furthermore, it is easy to check that \eqref{elas_S_0} is correct. Hence we proved this lemma.
\endpf

\begin{lemma}\label{elas_lemmaR} Consider  $C^1$ solutions to \eqref{LCEuler} with $c(x)$ given in $\eqref{elas1_c0}$ and prescribed initial data $\rho_0\equiv1$, $u_0=c(x)$. Assume $S(x,t)>0$ for any $(x,t)\in\mathbb R\times [0,T)$. There exist constants $\kappa_1$ and $\kappa_2$ depending on $\e$ such that
\begin{equation}
\label{elas_R_posi}
-\kappa_1<R(x,t)\leq 0,
\end{equation}
and
\begin{equation}
\label{elas_S_posi0}
0<\kappa_2<S(x,t),
\end{equation}
for any $(x,t)\in\mathbb R\times [0,T)$.

\end{lemma}
\pf
We first estimate $R$.
We have already proved \eqref{elas_R_posi} in $\Omega_{\{1\}}$, i.e. $\Omega_\xi$ with $\xi=1$,  in the left of $\Gamma$ in Lemma
\ref{elas_lemmaT}.
For $(x,t)$ in the right of forward characteristic $\Gamma_1$, $R=0$ since $c(x)$ is a constant when $x\geq x_{P_5}=1+\e^5$ and \eqref{cx0}.
For any smooth solutions, by \eqref{vac-eq18},
\[
R_t+S R_x \geq -\frac{c'}{4c} S^2.
\]
Using the same argument as in  \eqref{elas_R_est} and the initial energy
in the domain of dependence including the region between $\Gamma$ and $\Gamma_1$ is finite, we have 
\[
R>-\kappa_1.
\]
for some positive constant $\kappa_1$.\vspace{.1cm}

Now we proceed to prove \eqref{elas_S_posi0}. To the right (or left) of the vertical line passing $P_5$ (or $P_2$) which are both backward characteristics, $S$ equals to its initial constant data, hence \eqref{elas_S_posi0} is correct.
To consider $S$ in the region $\Omega_{\{1\}}$, i.e.  to the left of the characteristic $\Gamma$ starting from $P_4$, and to the right of the vertical line passing $P_2$,
by Lemmas \ref{lemma2.1} and Lemma \ref{elas_lemmaT}, on any backward characteristic $\Gamma_\xi$ when $t<T$ starting from the point $P_\xi(\xi,0)$,
\[
S_t+R S_x \geq \frac{c'}{4c^2}(S^2-a_{\xi}^2).
\]
with initial data $S_0(\xi)=\frac{2\e}{1-\e\xi}>a_\xi$.
Hence $S\geq S_0(\xi)$ on $\Gamma_\xi$ when $t<T$. So \eqref{elas_S_posi0} is correct  in this region for some $\kappa_2$.

On the backward characteristic starting from $[P_4,P_5]$, we also have that $S(x,t)$
has positive lower bound as the previous paragraph, since $S_0(\xi)\approx\frac{2\e}{1-\e}>a_{\{\xi=1\}}\approx a_\xi$ with differences at most in $O(\e^{5})$ by \eqref{elas_ineq_d} when $\xi\in(1,1+\e^{5})$. This completes the proof of \eqref{elas_S_posi0}, hence the proof of the lemma.
\endpf
\vspace{.2cm}

Now we proceed to find the blowup of $S$.
\begin{lemma}\label{elas_lemmaS}
Consider $C^1$ solutions to \eqref{LCEuler} with $c(x)$ given in $\eqref{elas1_c0}$ and prescribed initial data $\rho_0\equiv1$, $u_0=c(x)$. Assume $S(x,t)>0$ for any $(x,t)\in\mathbb R\times [0,T)$. 
 There exist positive constants $M_0$ such that 
\begin{equation}
\label{elas_S_posi}
0<\kappa_2<S(x,t)<+\infty,
\end{equation}
for any $(x,t)\in\mathbb R\times[0,T)$ with $T=M_0 \e^{-2}$ and 
$$
\lim\limits_{(x,t)\rightarrow (x^*,T)}\,S(x,t)=+\infty
$$
for some $x^*$ such that $0<1-x^*\leq O(1)\e^{-1}$.

\end{lemma}

\pf 
By Lemma \ref{lemma2.1}, we have
\beq\label{elass_ineq1}
 S_t+R S_x \leq \frac{d_2}{2}S^2.
\eeq
where $d_2=O(\e)$ and $0<S_0\leq O(\e)$. So $S$ stays finite until before 
$$\bar{t}=M_2\e^{-2},$$
for some constant $M_2>0$.

Then we show the blowup happens at a time in $O(\e^{-2})$.
For simplicity, we only consider the backward characteristic $\Gamma$ starting from
the $P_4(1,0)$ on $(x,t)$-plane.
For any $(x,t)$ on $\Gamma$, by the estimate of $R$ in Lemma \ref{elas_lemmaR}, we have
\beq\label{elas_xtest}
\frac{1-x}{t}\leq a,
\eeq
where we use
\[
a=\frac{\e}{2(1-\e)},\quad \text{to\ denote}\quad
a_\xi=\frac{\e}{2(1-\e\xi)}\quad \text{at}\quad \xi=1.
\]

Before $\Gamma$ interacts with $\Gamma_0$ which will happen not earlier than $O(\e^{-9})$ by \eqref{elas_xtest},
by Lemma \ref{lemma2.1}, Lemma \ref{elas_lemmaR}, \eqref{elas_xtest} and definition of $c(x)$, 
\beq\label{elas_s_ineq}
S_t+R S_x \geq \frac{c}{4}(S^2-a^2)\geq \frac{\e}{4(1+a\e t-\e)}(S^2-a^2).
\eeq
Studying the ODE
\beq\label{elas_ds}
\frac{dg}{dt^-}= \frac{\e}{4(1+a\e t-\e)}(g^2-a^2)
\eeq
with initial data 
\beq\label{elas_s00}
g(0)=S_0(1)=\frac{2\e}{1-\e}>a,
\eeq
one has that $g$ blows up at $t^*$, which satisfies
$$
2\ln\left(\frac{S_0(1)+a}{S_0(1)-a}\right)+\ln(1-\e)=\ln(1+a\e t^*-\e).
$$
Therefore
\beq\label{elas_tstar}
t^*=\frac{16}{9a\e}(1-\e)=\frac{32(1-\e)^2}{9\e^2}=M_3\e^{-2},
\eeq
for constant $M_3=\frac{32(1-\e)^2}{9}$. By comparison theorem of ODE, $S(x,t)$ blows up not later than $t^*=M_3\e^{-2}$. 

Therefore, there exists $M_0\in [M_2, M_3]$ such that for any $(x,t)\in\mathbb R\times[0,T)$ with $T=M_0\e^{-2}$
$$0<S(x,t)<+\infty,\quad
\lim\limits_{(x,t)\rightarrow (x^*,T)}\,S(x,t)=+\infty
$$
for some $x^*$. By \eqref{elas_xtest} and $T=M_0\e^{-2}$, we have
$0<1-x^*\leq O(1)\e^{-1}$.
Hence we complete the proof of the lemma by \eqref{elas_S_posi0}.
%
%

\endpf

Next we show that the assumption that $S>0$ is true for all $C^1$ solutions in our initial value problems, which implies that $S$ is uniformly bounded away from zero by Lemma \ref{elas_lemmaR}.

\begin{lemma}{\label{elas_Sgeq0}}
For any $C^1$ solutions to \eqref{LCEuler} with $c(x)$ given in $\eqref{elas1_c0}$ and prescribed initial data $\rho_0\equiv1$, $u_0=c(x)$, one has $S(x,t)>0$ for any $(x,t)\in\mathbb R\times [0,T)$. Hence Lemmas \ref{lemma0}$\sim$\ref{lemmaS} are correct without the assumption that $S(x,t)>0$ in the beginning. 

\end{lemma}
\pf
Note $R=0$ and $S$ are positive constants in the left of $\Gamma_0$ and in the right $\Gamma_1$ respectively, since $c(x)$ has constant value on each of these two regions. Denote the finite region between these two regions as $\Omega^*$ with $t<T$.

We prove the lemma by contradiction. Assume that $S=0$ somewhere. Then
$S=0$ must first happen in $\Omega^*$. We could find the minimum time such that $S=0$. Assume that $S(\hat x, \hat T)=0$ for some point $(\hat x, \hat T)$ in
$\Omega^*$ and $S(x,t)>0$
for any  $(x,t)\in\mathbb R\times[0,\hat T)$. Then running the proofs in Lemmas \ref{elas_lemma0}$\sim$\ref{elas_lemmaR}, we can still get \eqref{elas_S_posi} for $(x,t)\in\mathbb R\times[0,\hat T]$ which contradicts to $S(\hat X, \hat T)=0$. Hence, $S(x,t)>0$ for any $(x,t)\in\mathbb R\times[0,T)$. 
\endpf\vspace{.2cm}

Finally we are ready to prove Theorem \ref{main-elast}.
\vspace{.2cm}

{\noindent{\bf Proof of Theorem \ref{main-elast}.}}
Combining the a priori estimates in Lemma \ref{elas_lemmaR} and Lemma \ref{elas_lemmaS}, using \eqref{vac-eq4}, we know the a priori condition (A) is true for any $T^*<T=M_0\e^{-2}$,
where we use the fact: for any $C^1$ solution, $S$ is bounded above in the closed set $\Omega^*$ defined in the previous lemma with $t\in[0,T^*]$ since $S$ is not infinity when  $t\in[0,T^*]$, and $S$ has constant value in the left of $\Gamma_0$ or in the right of $\Gamma_1$.

As a conclusion, by Theorem \ref{ex-sm-sl}, the initial value problem of \eqref{LCEuler} with the prescribed initial data exists a unique $C^1$ solution $(\rho(x,t), u(x,t))$ when
$(x,t)\in\mathbb R\times[0,T)$. 
Furthermore, by Lemma \ref{elas_lemmaS}, $S$ is uniformly positive and
\[
\lim_{t\rightarrow{T^-}}S(x^*,t)=\infty
\]
for some $|x^*|=O(\e^{-1})$, while $R$ is  non-negative and uniformly bounded by Lemma \ref{elas_lemmaR}. 
So by \eqref{vac-eq4},
\beq
\lim_{t\rightarrow{T^-}}u(x^*,t)=\infty \quad \lim_{t\rightarrow{T^-}}\rho(x^*,t)=0 \quad  \lim_{t\rightarrow{T^-}}\rho u(x^*,t)=B\quad \text{and}\quad
\lim_{t\rightarrow{T^-}}E(x^*,t)=\infty,
\eeq
where $B$ is a finite constant.
Hence we complete the proof of Theorem \ref{main-elast}.

\endpf

\section{Singularity formations in duct flow for Chaplygin gas: proof of Theorem \ref{main} and \ref{main2}}
\setcounter{equation}{0}
In this section, without any ambiguity, we still use several notations used in the previous section, such as $\e$, $\Omega$, $M_i$ etcs.

To prove Theorem \ref{main}, we need first set up the initial data $\rho_0$ and $u_0$ 
corresponding to $c(x)$ in \eqref{c1}
\beq \label{rho0duct}
\rho_0(x)=\left\{
\begin{array}{rcl}
\e^{-\alpha-1},&&x\in(-\infty,-2\e),\\
\phi_1(x),& &x\in[-2\e,-\e),\\
1, & &x\in[-\e,\e],\\
\phi_2(x),& &x\in(\e,2\e],\\
\e^{-\alpha-1},&&x\in(2\e,\infty),\\
\end{array}\right.
\eeq
and
\beq \label{u0duct}
u_0(x)=\frac{c(x)}{\rho_0(x)} -\e^{\alpha+1},
\eeq
where the function $\phi_1$ and $\phi_2$ are $C^1$ increasing and decreasing functions
on $[-2\e,-\e)$ and $(\e,2\e]$ respectively.  We collect some useful information on initial data here. 
First $\rho_0$ and $c(x)$ are $C^2$ functions bounded away from zero and infinity. $|c'(x)|$ is bounded above. $u_0$ is a $C^1$ function with finite upper and lower bounds. For any $x$,
\beq\label{c15}
2\leq c(x)\leq 4,
\eeq

\beq\label{c_derivative}
 c'(x)=\left\{
\begin{array}{rcl}
0,&&x\in(-\infty,-1-\eta),\\
\psi'_1(x)\in[0,\e^\alpha],&& x\in[-1-\eta,-1),\\
\e^\alpha, & &x\in[-1,1],\\
\psi'_2(x)\in[0,\e^\alpha],&& x\in(1,1+\eta],\\
0,&&x\in(1+\eta,\infty),\\
\end{array}\right.
\eeq\vspace{.1cm}

\beq \label{1overrho}\frac{1}{\rho_0(x)}=\left\{
\begin{array}{rcl}
\e^{\alpha+1},&&x\in(-\infty,-2\e),\\
\frac{1}{\phi_1(x)}\in[\e^{\alpha+1},1],&& x\in[-2\e,-\e),\\
1, & &x\in[-\e,\e],\\
\frac{1}{\phi_2(x)}\in[\e^{\alpha+1},1],& &x\in(\e,2\e],\\
\e^{\alpha+1},&&x\in(2\e,\infty),\\
\end{array}\right.
\eeq
and
\beq\label{RS0}
R_0(x)=-\e^{\alpha+1}\quad\text{and}\quad \e^{\alpha+1}<S_0(x)<+\infty.
\eeq

The constants $\alpha\in[0,1)$ and $0<\eta<\e^3$. Furthermore $0<\e\ll1$ is a small given number which will be provided in the proof of the theorem.

\begin{remark} \label{remarkRho0}
(i) The construction of the initial data $\rho_0$, $u_0$ and $c(x)$ is also motivated by the seminal work \cite{ghz}. The initial data given by equations (1.4) and (1.5) in Theorem 1 of \cite{ghz} were constructed for their unknown state $\phi$ (which is equivalent to $x(X,t)$ in our equation \eqref{LCEuler} if $\rho_0=1$). 
\\
(ii) When $\rho_0\not\equiv 1$, \eqref{LCEuler} is not equivalent to \eqref{sn6} any more. This requires us to do extra constructions on $c$ (as a function of $x\in(-\infty,\infty)$). While in \cite{ghz}, the wave speed $c$ was a function of their unknown state $\phi$ (see (1.1) in \cite{ghz}).  They assumed $c$ to be a uniformly positive and bounded smooth function.
\\
(iii)The constructions on $R_0$, $S_0$ in \eqref{RS0} play an essential role in our proof.   
\end{remark}

\begin{figure}[htb]
\centering
\includegraphics[width=12cm,height=6cm]{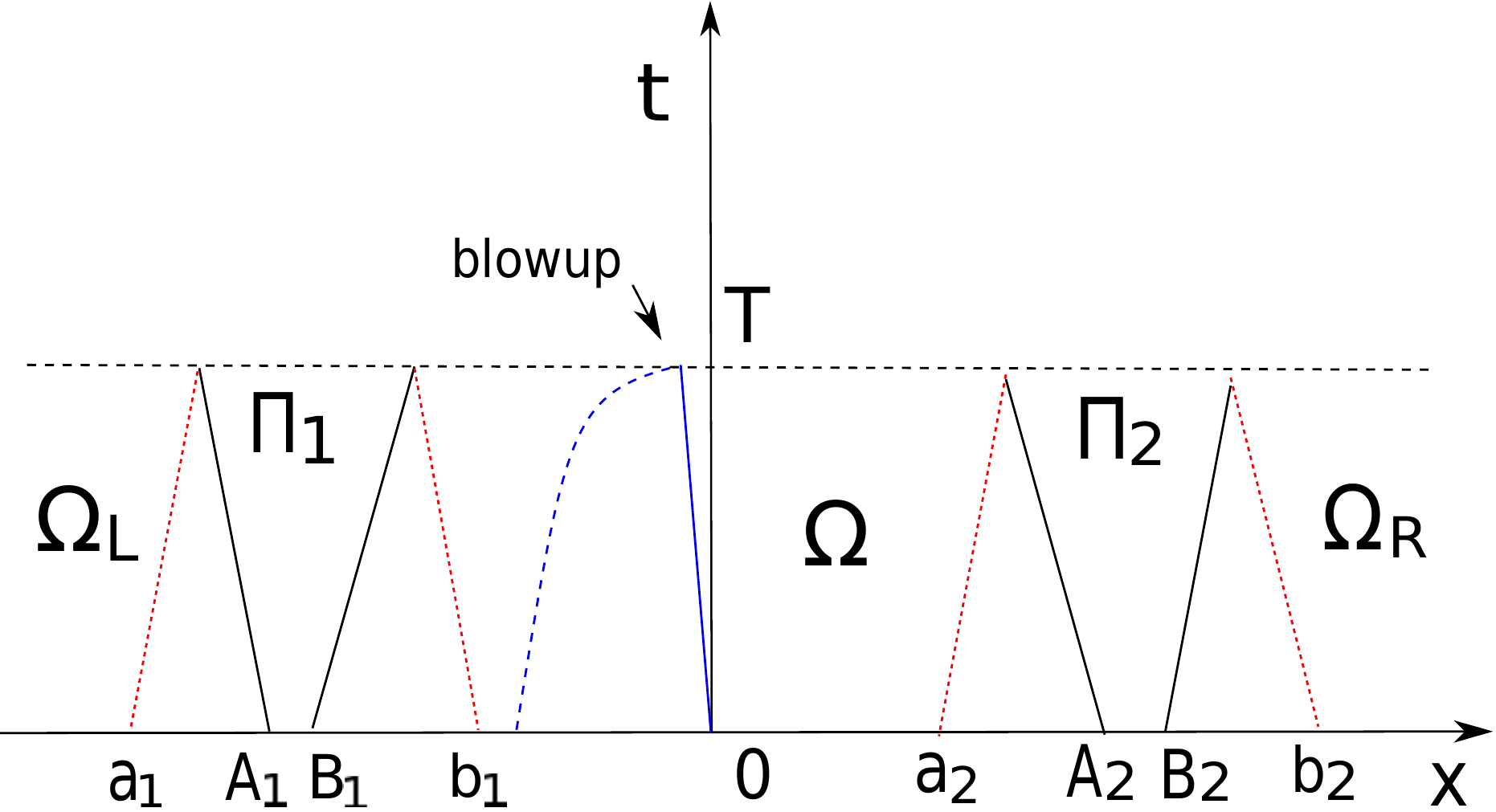}
\caption{Proof of Theorem \ref{main}}
\label{euler}
\end{figure}

Assuming that $\rho(x,t),~u(x,t)$ are $C^1$ solutions to Cauchy problem of \eqref{LCEuler} with given $c$, $\rho_0$ and $u_0$, we first do some analysis on the domains of dependence of different pieces of the initial data (see Figure \ref{euler}), where $A_1=-1-\eta$, $B_1=-1$,
 $A_2=1$, $B_2=1+\eta$ and when $ t\in [0, T^*]$ for any $T^*<T$ with
 \beq\label{T_def}
 T=M_1\e^{-\alpha}, 
 \eeq
 where $M_1$ is a fixed constant that will be provided later. In fact, $T$ is the time of blowup.  
 
We use $\Omega$,
$\Omega_L$, $\Omega_R$ to denote three domains of dependence respect to $[-1,1]$, $(-\infty,-1-\eta)$, $(1+\eta,+\infty)$ respectively (boundaries are black solid lines). $\Pi1$ and $\Pi2$
are two regions in between
$\Omega_L$,  $\Omega$ and $\Omega_R$. The domains of dependence $\Omega_{\Pi1}$ and $\Omega_{\Pi2}$
including $\Pi1$ and $\Pi2$ respectively are regions with red dash lines as boundaries (with initial bases $[a_1,b_1]$ and $[a_2,b_2]$ respectively).

We first have a lemma for $C^1$ solutions.

\begin{lemma}\label{lemma0}
For any $C^1$ solutions with  prescribed initial data,
 \[
 R<0.
 \]
\end{lemma}
\pf
By  Lemma \ref{lemma2.1}, for any $x\in(-\infty,\infty)$
\beq\label{Rin0}
R_t+S R_x \leq\frac{1}{4} R^2,
\eeq
where left hand side is derivative along a minus characteristic.
By \eqref{RS0} and ODE comparison theorem, 
\[R<0\]
 for any $C^1$ smooth solution.
\endpf
\vspace{.2cm}

In Lemmas \ref{lemmaT}$\sim$\ref{lemmaS}, we restrict our consideration on the $C^1$ solutions with $S(x,t)>0$ for any $(x,t)\in\mathbb R\times [0,T)$ where $T$ is defined in \eqref{T_def}. For these solutions, the plus and minus characteristic go in forward and backward directions respectively. 
So we can always use Figure \ref{f0} to study the propagation of the solution along characteristics. 
We will show that this a priori condition on $S$ is satisfied for any $C^1$ solutions in our initial value problems in Lemma \ref{Sgeq0}.

\begin{lemma}\label{lemmaT}
Consider  $C^1$ solutions with  prescribed  initial data and $S(x,t)>0$ for any $(x,t)\in\mathbb R\times [0,T)$.
Let $T_0$ be the intersection time of forward characteristic from $(-1,0)$ and backward characteristic from $(1,0)$. Then there exists a constant $K_1>0$ such that 
\begin{equation}
\label{est-T0}
T_0\geq K_1\e^{-1}>T.
\end{equation}
And we also have the estimates
\begin{equation}
\label{distBb}
|b_1-B_1|+|A_2-a_2|\leq K_2\e^{\frac{1-\alpha}{2}},
\end{equation}
where $b_1$, $B_1$, $A_2$ and $a_2$ are on Figure \ref{euler}.

\end{lemma}

\pf
It can be calculate by \eqref{vac-eq4}, \eqref{c1}, \eqref{u0duct} and \eqref{1overrho} that the initial energy in $\Omega$ satisfies
\beq\label{Ef}
\int_{-1}^{1}\rho_0(S^2+R^2)(x,0)\,dx
=2\int_{-1}^{1}\rho_0^{-1}(\rho_0^2u_0^2+c^2)\,dx
=2\int_{-1}^{1}\rho_0^{-1}(\left(c-\rho_0\e^{\alpha+1}\right)^2+c^2)\,dx
\leq K_3\e.
\eeq
For any $x_1, x_2\in[-1,1]$, let $(x_0,t_0)$ be the intersection of forward and backward characteristics $\gamma_{+}$, $\gamma_{-}$ starting $(x_1,0)$ and $(x_2,0)$ respectively.
By Lemma \ref{vac-lemma1} of finite propagation and \eqref{Ef}, we know 
%
\begin{equation}
\label{aa1}
\begin{split}
&2\int_{0}^{t_0}cS^2(x_{+}(t),t)\,dt+2\int_{0}^{t_0}cR^2(x_{-}(t),t)\,dt\\
=&2\int_{x_1}^{x_0}cS(x,t_{+}(x))\,dx-2\int_{x_0}^{x_2}cR(x,t_{-}(x))\,dx\leq K_3 \e.
\end{split}
\end{equation}
where $K_i$ and $M_i$ always mean positive constant in this paper.

Therefore
\begin{equation}
\label{aa2}
\begin{split}
&x_2-x_1=x_0-x_1+x_2-x_0\\
=& \int_{0}^{t_0}S(x_{+}(t),t)\,dt-\int_{0}^{t_0}R(x_{-}(t),t)\,dt\\
\leq& \left(\int_{0}^{t_0} \frac{1}{c}dt\right)^\frac{1}{2} \left(2{\int_{0}^{t_0}cS^2(x,x_{+}(t))\,dt}\right)^\frac{1}{2}+\left(\int_{0}^{t_0} \frac{1}{c}dt\right)^\frac{1}{2} \left(2{\int_{0}^{t_0}cR^2(x,x_{-}(t))\,dt}\right)^\frac{1}{2}\\
\leq& K_4 t_0 ^{\frac{1}{2}} \e^{\frac{1}{2}},
\end{split}
\end{equation}
It is easy to see that \eqref{est-T0} by \eqref{aa2}.

Repeat the proof of \eqref{aa2}, we can prove
\begin{equation}
|b_1-B_1|+|A_2-a_2|\leq 2K_5 T ^{\frac{1}{2}} \e^{\frac{1}{2}}\leq K_2\e^{\frac{1-\alpha}{2}},
\end{equation}
where in fact set $(x_1,t_1)=B_1$ or $a_2$ in Figure \ref{f0} then when  
$x_0-x_1=O(\e^{\frac{1-\alpha}{2}})$,
 \begin{eqnarray}
x_0-x_1&=& \int_{t_1}^{t_0}S(x_{+}(t),t)\,dt\nn\\
		&\leq&\frac{1}{\sqrt{2}} \left(\int_{t_1}^{t_0} \frac{1}{c}dt\right)^\frac{1}{2} \left({\int_{t_1}^{t_0}2cS^2(x,x_{+}(t))\,dt}\right)^\frac{1}{2}\nn\\		
		&=& \frac{1}{\sqrt{2}}\left(\int_{t_1}^{t_0} \frac{1}{c}dt\right)^\frac{1}{2} \left(2\int_{x_1}^{x_0}cS(x,t_{+}(x))\,dx\right)^\frac{1}{2}\nn\\
		&\leq& \frac{1}{\sqrt{2}}\left(\int_{t_1}^{t_0} \frac{1}{c}dt\right)^\frac{1}{2} \left(\int_{x_1}^{ x_2}\rho_0(S^2_0+R^2_0)\,dx\right)^\frac{1}{2}\nn\\
		&\leq& M_4\cdot (t_0 - t_1)^{\frac{1}{2}}(x_0-x_1)^{\frac{1}{2}} \e^{\frac{1+\alpha}{2}}\nn
\end{eqnarray}
so $t_0 - t_1\geq O(\e^{-\alpha})$.

Hence we complete the proof of Lemma \ref{lemmaT}. 
\endpf
\vspace{.2cm}

\begin{remark}
\label{remark4.3}
(i) By \eqref{distBb}, it is easy to get that the initial energy in $\Omega_{\Pi_1}$ or $\Omega_{\Pi_2}$ is not greater than $O(\e^{2+\alpha})$ by \eqref{distBb}.\\
(ii)  On regions to the left or right of $\Omega$, $\Omega_{\Pi_1}$ and $\Omega_{\Pi_2}$, the wave speeds $S$ and $R$ are constants by \eqref{cx0} ($c(x)$ has constant values in those regions) and $S$ and $R$ are constants initially
in each of these two regions.
\end{remark}

\begin{lemma}\label{lemmaR} Consider  $C^1$ solutions with prescribed initial data and $S(x,t)>0$ for any $(x,t)\in\mathbb R\times [0,T)$. There exist constants  $\delta_2$ and $\delta_3$ such that
\begin{equation}
\label{R_posi}
-\delta_3<R(x,t)<-\delta_2<0.
\end{equation}

\end{lemma}

\pf By Remark \ref{remark4.3} (ii) and \eqref{RS0}, we only need to consider the solution on $\Omega$, $\Omega_{\Pi_1}$ and $\Omega_{\Pi_2}$.
By  Lemma \ref{lemma2.1}, for any $x\in(-\infty,\infty)$
\beq\label{Rin}
R_t+S R_x \leq\frac{1}{8} R^2.
\eeq
By \eqref{RS0} and ODE comparison theorem, there exists $\delta_2$ such that 
\[R<-\delta_2<0.\]
For any forward characteristic in $\Omega\cup\Omega_{\Pi_1}\cup\Omega_{\Pi_2}$, by  \eqref{vac-eq12}, \eqref{c1} and \eqref{c_derivative} we have
\beq\label{R_ineq}
R_t+S R_x \geq - \frac{\e^\alpha}{16}cS^2.
\eeq
Integrating it along any forward characteristic $\gamma_+$, by  Lemma \ref{vac-eq18}, Remark \ref{remark4.3} (i) and \eqref{RS0}, we have
\begin{eqnarray}
R(t_0,x_0)&\geq& - \frac{\e^\alpha}{16}\int_0^{t_0} cS^2(t,x_{+}(t))dt +R_0( x_1)\nn\\
&=&- \frac{\e^\alpha}{16}\int_{x_1}^{x_0} cS(t_{+}(x),x)dx-\e^{\alpha+1}\nn\\
&\geq &- \frac{\e^\alpha}{16}\int_{x_1}^{x_2} \rho_0(S_0^2+R_0^2)dx-\e^{\alpha+1}\nn\\
&\geq& -K_6 \e^{1+\alpha}:=-\delta_3.\label{R_est}
\end{eqnarray}
\endpf
\vspace{.2cm}

Now we proceed to estimate $S$.
\begin{lemma}\label{lemmaS}
Consider  $C^1$ solutions with  prescribed  initial data and $S(x,t)>0$ for any $(x,t)\in\mathbb R\times [0,T)$. 
 There exist positive constants $M_1$ and $\delta_1$ such that 
\begin{equation}
\label{S_posi}
0<\delta_1<S(x,t)<+\infty,
\end{equation}
for any $(x,t)\in\mathbb R\times[0,T)$ with $T=M_1 \e^{-\alpha}$ and 
$$
\lim\limits_{(x,t)\rightarrow (x^*,T)}\,S(x,t)=+\infty
$$
for some $x^*=O(1)\e^{\frac{1-\alpha}{2}}$.

\end{lemma}

\pf By Remark \ref{remark4.3} (ii) and \eqref{RS0}, we only need to consider the solution on $\Omega$, $\Omega_{\Pi_1}$ and $\Omega_{\Pi_2}$. 

By Lemma \ref{lemma2.1} and the estimate of $R$ in Lemma \ref{lemmaR}, we have
\beq\label{s_ineq1}
 S_t+R S_x \leq \frac{\e^{\alpha}}{8}S^2.
\eeq
So $S$ stays finite until before 
$$\bar{t}=M_2\e^{-\alpha},$$
where $M_2=\frac{8}{S_0(x)}$.

When $x\not\in[-2\e,2\e]$, we have $S_0(x)=O(\e^{1+\alpha})$. So by the comparison theorem of ODE, $S(x,t)$ along any backward characteristic starting from an initial point with $x\not\in[-2\e,2\e]$, will not blowup until $\bar{t}=O(\e^{-1-2\alpha})$. Comparing to the blowup time in $O(\e^{-\alpha})$ proved later, we see that the blowup can only happen on a backward characteristic starting from the initial interval $[-2\e,2\e]$. 

Then we show the blowup happens at a time in $O(\e^{-\alpha})$.
For simplicity, we only consider the backward characteristic $\Gamma$ starting from
the origin on $(x,t)$-plane.
On $\Gamma$, by Lemma \ref{lemma2.1} and the estimate of $R$ in Lemma \ref{lemmaR}, 
\beq\label{s_ineq}
S_t+R S_x \geq \frac{\e^\alpha}{16}S^2-\frac{\e^\alpha}{8}\delta_3^2.
\eeq
Studying the ODE
\beq\label{ds}
\frac{dg}{dt^-}= \frac{\e^\alpha}{16} g^2-\frac{\e^\alpha}{8}\delta_3^2
\eeq
with initial data 
\beq\label{s00}
g(0)=S_0(0)=6-\e^{\alpha+1}\in(5,6),
\eeq
one has that $g$ blows up at
\beq\label{tstar}
t^*=\frac{8}{K_7\e^{1+2\alpha}}\, \ln\left| \frac{S_0(0)+2\delta_3 } {S_0(0)-2\delta_3} \right|= M_3\e^{-\alpha},
\eeq
for some constant $M_3>0$. By comparison theorem of ODE, $S(x,t)$ blows up not later than $\bar{t}=M_3\e^{-\alpha}$. 

Therefore, there exists $M_1\in [M_2, M_3]$ such that for any $(x,t)\in\mathbb R\times[0,T)$ with $T=M_1\e^{-\alpha}$
$$0<S(x,t)<+\infty,\quad
\lim\limits_{(x,t)\rightarrow (x^*,T)}\,S(x,t)=+\infty
$$
for some $x^*$. By \eqref{aa2} and $T=M_1\e^{-\alpha}$, we have
$x^*=O(1)\e^{\frac{1-\alpha}{2}}$.

For any $(x,t)\in\mathbb R\times[0,T)$, by \eqref{s_ineq}, we have
\[
S_t+R S_x \geq -\frac{\e^\alpha}{8}\delta_3^2,
\]
then by \eqref{RS0} ($S_0(x)>\e^{\alpha+1}$) and comparison theorem of ODE, 
we obtain 
\beq\label{S_posi1}
0<\delta_1<S(x,t)<+\infty
\eeq
for some constant $\delta_1$, which completes the proof of the lemma.

\endpf

Next we show that the assumption that $S>0$ is true for all $C^1$ solutions in our initial value problems.

\begin{lemma}{\label{Sgeq0}}
For any $C^1$ solutions with  prescribed  initial data,  $S(x,t)>0$ for any $(x,t)\in\mathbb R\times [0,T)$. Hence Lemmas \ref{lemma0}$\sim$\ref{lemmaS} are correct without the assumption that $S(x,t)>0$ in the beginning.
\end{lemma}
\pf
We prove it by contradiction. Assume that $S=0$ somewhere.

Note $R$ and $S$ are non-zero negative and positive constants respectively when $(x,t)\not\in\Omega\cup\Omega_{\Pi_1}\cup\Omega_{\Pi_2}$. So if $S=0$ then
it must first happen in $\Omega\cup\Omega_{\Pi_1}\cup\Omega_{\Pi_2}$. We could find the minimum time such that $S=0$. Assume that $S(\hat x, \hat T)=0$ for some
$(\hat x, \hat T)\in\Omega\cup\Omega_{\Pi_1}\cup\Omega_{\Pi_2}$ and $S(x,t)>0$
for any  $(x,t)\in\mathbb R\times[0,\hat T)$. Then running the proofs in Lemmas \ref{lemma0}$\sim$\ref{lemmaS}, we can still get \eqref{S_posi1} for $(x,t)\in\mathbb R\times[0,\hat T]$ which contradicts to $S(\hat X, \hat T)=0$. Hence, $S(x,t)>0$ for any $(x,t)\in\mathbb R\times[0,T)$.
\endpf\vspace{.2cm}

Finally we are ready to prove Theorem \ref{main}.
\vspace{.2cm}

{\noindent{\bf Proof of Theorem \ref{main}.}
Combining the a priori estimates in Lemma \ref{lemmaR} and Lemma \ref{lemmaS}, using \eqref{vac-eq4}, we know the a priori condition (A) is true for any $T^*<T=M_1\e^{-\alpha}$,
where we use the fact: for any $C^1$ solution, $S$ is bounded above in the closed set $\Omega\cup\Omega_{\Pi_1}\cup\Omega_{\Pi_2}$ with $t\in[0,T^*]$ since $S$ is not infinity when  $t\in[0,T^*]$, then uniformly bounded from above for any $(x,t)\in\mathbb R\times[0,T^*]$ by Remark \ref{remark4.3} (ii).

As a conclusion, by Theorem \ref{ex-sm-sl}, there exists a unique $C^1$ solution $(\rho(x,t), u(x,t))$ for equation \eqref{LCEuler} when
$(x,t)\in\mathbb R\times[0,T)$ with the prescribed initial data. 

Furthermore, by Lemma \ref{lemmaS}
\[
\lim_{t\rightarrow{T^-}}S(x^*,t)=\infty
\]
while $R$ is uniformly bounded and negative by Lemma \ref{lemmaR}. So by \eqref{vac-eq4},
\beq
\lim_{t\rightarrow{T^-}}u(x^*,t)=\infty \quad \lim_{t\rightarrow{T^-}}\rho(x^*,t)=0 \quad  \lim_{t\rightarrow{T^-}}\rho u(x^*,t)=B\quad \text{and}\quad
\lim_{t\rightarrow{T^-}}E(x^*,t)=\infty,
\eeq
where $B$ is a finite constant.
Here $|x^*|=O(\e^{\frac{1-\alpha}{2}})<O(1)$ and the blowup must be on some characteristic starting form the initial interval $x\in[-2\e, 2\e]$. Hence we complete the proof of Theorem \ref{main}.

\endpf
%

\bigskip
Finally we prove the Theorem \ref{main2} for the
Cauchy problem of  \eqref{3deuler} with radially symmetry, i.e. \eqref{LCEuler'} or  \eqref{LCEuler}  
with $c(x)=x^m$ ($m=1$, cylindrical symmetric solution; $m=2$, spherically symmetric solution) and the radius $x\in[1,3]$. This theorem is also correct when $x\in[a,b]$ with $0<a<b$.

\bigskip
\noindent{\bf Proof of Theorem \ref{main2}}.} When $x\in[1,3]$, we consider initial data $\rho_0(x)$ and $u_0(x)$ satisfying
\begin{itemize}

\item[(1).] 
\beq \rho_0(x)=\left\{
\begin{array}{rcl}
\e^{-1},&&x\in[1,2-2\e),\\
\phi_3(x),&& x\in[2-2\e,2-\e),\\
1, && x\in[2-\e,2+\e],\\
\phi_4(x),&& x\in(2+\e,2+2\e],\\
\e^{-1},&&x\in(2+2\e,3],\\
\end{array}\right.
\eeq
\item[(2).] $u_0(x)=\frac{c(x)}{\rho_0(x)} -\e.$
\end{itemize}
The function $\phi_3$ and $\phi_4$ are increasing and decreasing smooth functions
on $[2-2\e,2-\e)$ and $(2+\e,2+2\e]$ respectively. Furthermore $\e<1$ is a small given number which will be provided in the proof of the theorem.

In order to directly use the proof for Theorem \ref{main}, we extend the definition of initial data from $x\in[1,3]$ to $x\in(-\infty,\infty)$
by
\begin{itemize}
\item[(1').] 
\beq
\label{ctilde}
 \tilde c(x)=\left\{
\begin{array}{rcl}
1-\delta,&&x\in(-\infty,1-\eta),\\
\psi_3(x),& &x\in[1-\eta,1),\\
x^m, & &x\in[1,3],\\
\psi_4(x),& &x\in(3,3+\eta],\\
3^m+\delta,&&x\in(3+\eta,\infty),\\
\end{array}\right.
\eeq
where $\psi_3(x)$ is an increasing $C^2$ convex positive function on $x\in[1-\eta,1)$ and $\psi_4(x)$ is an increasing $C^2$ concave positive function on $x\in(3,3+\eta]$. The positive constants 
\[\delta\ll\eta\ll\e\ll1.\] 
\item[(2').] 
\beq \tilde\rho_0(x)=\left\{
\begin{array}{rcl}
\e^{-1},&&x\in(-\infty,2-2\e),\\
\phi_3(x),&& x\in[2-2\e,2-\e),\\
1, & &x\in[2-\e,2+\e],\\
\phi_4(x),&& x\in(2+\e,2+2\e],\\
\e^{-1},&&x\in(2+2\e,\infty),\\
\end{array}\right.
\eeq
\item[(3').] $\tilde u_0(x)=\frac{\tilde c(x)}{\tilde \rho_0(x)} -\e.$
\end{itemize}

The initial data in Theorem \ref{main2} are very similar to the initial data in Theorem \ref{main} with $\alpha=0$. 
In fact, we only need to change $x\in[-1,1]$ to radius $x\in[1,3]$ and slightly change the values of $c$,
then we can prove Theorem \ref{main2} by an entirely same way as the proof in Theorem \ref{main},
and finally we only have to use the piece of solution on $\Omega_{symm}$ with $t\in[0,T)$ after finding the $C^1$ solution for $(x,t)\in \mathbb R \times [0,T)$. We leave the details to the reader.

\bigskip

\section*{Acknowledgement} 
The work of Huang and Liu is partially supported by NSF-DMS grants 
1216938 and 1109107. The authors would like to thank Professors Alberto Bressan, Heldge Kristian Jenssen, Jiequan Li and Yuxi Zheng for their helpful comments and suggestions.

%


\begin{thebibliography}{99}
\bibitem{arnold}
V.~I.~Arnolʹd, Mathematical methods of classical mechanics. Translated from the Russian by K. Vogtmann and A. Weinstein. Second edition. Graduate Texts in Mathematics, 60. Springer-Verlag, New York, 1989.

\bibitem{AH} G.~Al\`{i} and J.~Hunter, Orientation waves in a director field
with rotational inertia, {\it Kinet. Relat. Models}, {\bf 2} (2009),
 1--37.

\bibitem{bj}P.~Baiti and H.~K.~Jenssen,
Blowup in $L\sp \infty$ for a class of genuinely nonlinear
   hyperbolic systems of conservation laws,
   {\it Discrete Contin. Dynam. Systems},
   {\bf 7}:4 (2001), 837--853.

\bibitem{Bressan}A.~Bressan, Hyperbolic Systems of Conservation
  laws: The One-dimensional Cauchy Problem, Oxford Lecture
  Ser. Math. Appl. 20, Oxford Univ. Press, Oxford, 2000.


\bibitem{BZ}A.  Bressan and Yuxi Zheng,  
Conservative solutions to a nonlinear variational wave equation,
\emph{Comm. Math. Phys.}, \textbf{266} (2006), 471--497.

\bibitem{G3}  G.~Chen,
 Formation of singularity and smooth wave propagation for the
non-isentropic compressible Euler
equations, {\em J. Hyperbolic Differ. Equ.}, \textbf{8}:4 (2011), 671--690.


\bibitem{G5} G.~Chen and R.~Young, Smooth waves and gradient blowup for the inhomogeneous 
wave equations, \emph{J. Differential Equations}, {\bf 252}:3 (2012), 2580--2595.

\bibitem{G6} G.~Chen and R.~Young,  Shock-free Solutions of the Compressible Euler 
Equations, {\em submitted}, available at http://arxiv.org/pdf/1204.0460v1.pdf


\bibitem{G8} G.~Chen, R.~Young and Q.~Zhang, {\it Shock formation in the compressible 
Euler equations and related systems}, J. Hyperbolic Differ. Equ., \textbf{10}:1 (2013), 149--172.

\bibitem{CZZ12} G.~Chen, P.~Zhang and Y.~Zheng, Energy Conservative Solutions to a Nonlinear Wave
System of Nematic Liquid Crystals,  {\it Comm. Pure Appl. Anal.}, \textbf{12}:3 (2013), 1445--1468.

\bibitem{GCZ} G.~Chen and Y.~Zheng, Existence and singularity to a wave system of 
nematic liquid crystals, {\em J. Math. Anal. Appl.}, \textbf{398} (2013), 170--188.

\bibitem{GG} G.-Q.~Chen and J.~Glimm, Global solutions to the compressible Euler 
equations with geometrical structure,  {\em Comm. Math. Phys.}, \textbf{180}:1 (1996), 153--193.


\bibitem{chap}S.~Chaplygin, On gas jets, {\it Sci. Mem . Moscow Univ. Math. Phys.,} 
{\bf 21} (1904), 1--121.

\bibitem{[6]} D.~Christodoulou and A.~Tahvildar-Zadeh,  On the regularity of
spherically symmetric wave maps, {\it Comm. Pure Appl. Math.}, \textbf{46} (1993), 1041--1091.


\bibitem{Dafermos} C.~M.~Dafermos, \emph{Hyperbolic Conservations laws in
Continuum Physics (third edition)}, Springer-Verlag, Heidelberg, 2010.

\bibitem{DCL} Y.~Du, G.~Chen and J.~Liu,
The almost global existence for a 3-D wave equation of nematic liquid-crystals, {\em submitted}. 

\bibitem{[10]} J.~L.~Ericksen and D.~Kinderlehrer (eds.),
Theory and Application of Liquid Crystals, IMA Volumes in
Mathematics and its Applications, {\em Vol.~5}, Springer-Verlag, New York
(1987).

\bibitem{ghz}  R.~Glassey, J.~Hunter and Y.~Zheng,
Singularities of a variational wave equation, {\it J. Differential Equations}, 
{\bf 129} (1996), 49--78.


\bibitem{HR} Helge, Holden and Xavier, Raynaud, Global semigroup of conservative solutions of the nonlinear, {\it Arch. Ration. Mech. Anal.}, \textbf{201}:3 (2011), 871--964

\bibitem{huang-liu-weber} T.~Huang, C.~Liu and F.~Weber, Eulerian description of variational wave equation. In preparation.


\bibitem{je} H.~K.~Jenssen, Blowup for systems of conservation laws, {\it SIAM J. Math. Anal.},
\textbf{31}:4 (2000), 894--908.

\bibitem{jy} H.~K.~Jenssen and R.~Young, Gradient driven and singular flux blowup of smooth
 solutions to
   hyperbolic systems of conservation laws, {\it J. Hyperbolic Differ. Equ.},
   \textbf{1}:4 (2004), 627--641.
   
   \bibitem{john}  F. John, { Formation of Singularities
in One-Dimensional Nonlinear Wave Propagation,}  {\it Comm. Pure Appl.
Math., }  \textbf{27} (1974), 377--405.

\bibitem{john0} F. John, { Formation of singularities
in elastic waves}, {\it Lecture Notes in Physics},   \textbf{195} (1984), 194--210.
   
   \bibitem{lax0} P. D.~Lax, Hyperbolic systems of conservation
laws, II,
{\it Comm. Pure Appl. Math.}, \textbf{10} (1957), 537-566.

\bibitem{lax}
P. D.~Lax, Development of singularities of solutions of nonlinear hyperbolic partial differential 
equations. {\it J. Math. Phys.}, \textbf{5} (1964), 611-613.

 \bibitem{lm} P.~G.~LeFloch and M.~Westdickenberg, 
 Finite energy solutions to the isentropic Euler equations
with geometric effects, {\it J. Math. Pures Appl.}, \textbf{88} (2007), 389--429.

\bibitem{lei-liu-zhou}
Z.~Lei, C.~Liu and Y.~Zhou, Global existence for a 2D incompressible viscoelastic model with small strain, {\em Commun. Math. Sci.}, \textbf{5}:3 (2007), 595--616.

\bibitem{li1}T.~T.~Li, \emph{Global classical solutions for quasilinear hyperbolic systems,}
RAM: Research in Applied Mathematics, 32, Masson, Paris, 1994.

\bibitem{Liyu} T.~T.~Li and W.~C.~Yu,
  Boundary value problems for quasilinear hyperbolic systems,
   Duke University Mathematics Series, V,
   Duke University Mathematics Department,
   Durham, NC, 1985.

\bibitem{lin-liu-zhang} F.~H.~Lin, C.~Liu and P.~Zhang, 
 On Hydrodynamics of Viscoelastic Fluids, {\em Comm. Pure Appl. Math.}, \textbf{58}:11 (2005), 1437--1471.

\bibitem{lin} L.~Lin,
On the vacuum state for the equations of isentropic gas dynamics,
{\it J. Math. Anal. Appl.}, \textbf{121}:2 (1987), 406--425. 

\bibitem{liu-walkington}

C.~Liu, N.~J.~Walkington, An Eulerian description of fluids containing visco-elastic particles, {\em Arch. Ration. Mech. Anal.} \textbf{159}:3 (2001), 229--252. 

\bibitem{liu3} T.-P. Liu, Transonic gas flow in a duct of varying area, {\it Arch. Ration. Mech. Anal.}, \textbf{80} (1982), 1--18

\bibitem{serre}
D.~Serre,
Multidimensional Shock Interaction
for a Chaplygin Gas,
{\it Arch. Rational Mech. Anal.}, \textbf{191} (2009), 539--577.

\bibitem{[40]} J.~Shatah, Weak solutions and development of singularities
in the $SU(2)$ $\sigma$-model, {\it Comm. Pure Appl. Math.}, \textbf{41} (1988), 459--469.

\bibitem{[41]} J.~Shatah and A.~Tahvildar-Zadeh,  Regularity of harmonic
maps from Minkowski space into rotationally symmetric manifolds,
{\it Comm. Pure Appl. Math.}, \textbf{45} (1992), 947--971.

\bibitem{yng} R.~Young, Blowup in hyperbolic conservation laws. 
{\it Contemp. Math.,} \textbf{327} (2003), 379--387.

\bibitem{ys}
R.~Young and W.~Szeliga, Blowup with small BV data in hyperbolic conservation laws. 
{\it Arch. Ration. Mech. Anal.,} {\bf 179} (2006) 31--54.

\bibitem{ZZ10} P.~Zhang and Y.~Zheng, Conservative
solutions to a  system of variational wave equations of nematic
liquid crystals, {\it Arch. Ration. Mech. Anal.}, \textbf{195} (2010),701--727.

\bibitem{ZZ11} P.~Zhang and Y.~Zheng, Energy Conservative Solutions to a One-Dimensional 
Full Variational Wave
System,  {\it Comm. Pure Appl. Math.}, \textbf{65}:5 (2012), 683--726.





\end{thebibliography}
\end{document}